\documentclass[10pt]{article}

\usepackage{graphics}
\usepackage{graphicx}

\usepackage[a4paper, left=35mm,right=35mm,top=34mm,bottom=34mm]{geometry}
\usepackage[utf8]{inputenc}
\usepackage[T1]{fontenc}
\usepackage[english]{babel}

\usepackage{enumerate}
\usepackage{graphicx}
\usepackage{hyperref}
\hypersetup{
    colorlinks=true,
    linkcolor=blue,
    filecolor=magenta,      
    urlcolor=cyan,
}
\usepackage{listings}
\usepackage{color}

\definecolor{dkgreen}{rgb}{0,0.6,0}
\definecolor{gray}{rgb}{0.5,0.5,0.5}
\definecolor{mauve}{rgb}{0.58,0,0.82}
\lstdefinelanguage{MRGC++}{%
  language=C++,
  morekeywords={T, U, MPI_Irecv, MPI_Isend, MPI_Allreduce, MPI_Waitall, Compute, Map, abs, max, Swap, MPI_Recv_init, MPI_Send_init, MPI_Startall, Copy, Init, InitRecv, InitSend, InitAllReduce, Send, Recv, AllReduce, Finalize, InitSnapshot, Snapshot, SwitchAsync, SnapReduce, MPI_Test, MPI_Start}
}
\usepackage{mathtools,amsthm,amssymb,amsfonts}
\usepackage{algorithm}
\usepackage{algorithmic}
\makeatother
\theoremstyle{plain}
\newtheorem{theorem}{Theorem}
\newtheorem{corollary}{Corollary}

\newtheorem{assumption}{Assumption}
\newtheorem{lemma}{Lemma}
\theoremstyle{definition}
\newtheorem{definition}[theorem]{Definition}
\theoremstyle{definition}
\newtheorem{notation}[theorem]{Notation}
\theoremstyle{remark}
\newtheorem*{remark}{Remark}

\usepackage{caption} 
\usepackage{fancyhdr}

\author{
  {\normalsize Guillaume Gbikpi-Benissan}\thanks{Academy of Engineering, RUDN University, Moscow, Russian Federation
    (correspondence, guibenissan@gmail.com).}
  \and
  {\normalsize Fr\'ed\'eric Magoul\`es}\thanks{Universit\'e Paris-Saclay, CentraleSup\'elec, Gif-sur-Yvette, France
    (correspondence, frederic.magoules@hotmail.com).}
}
\title{Resilient asynchronous primal Schur method}
\date{}

\begin{document}
\maketitle
\thispagestyle{fancy}

\begin{abstract}
\noindent This paper introduces the application of the asynchronous iterations theory within the framework of the primal Schur domain decomposition method. A suitable relaxation scheme is designed, which asynchronous convergence is established under classical spectral radius conditions. For the usual case where the local Schur complement matrices are not constructed, suitable splittings only based on explicitly generated matrices are provided. Numerical experiments are conducted on a supercomputer
for both Poisson's and linear elasticity problems.
The asynchronous Schur solver outperformed the classical conjugate-gradient-based one in case of compute node failures.
\end{abstract}


\section{Introduction}\label{sec1}

Asynchronous iterative methods are gaining more and more attention in the scientific computing field, especially for numerical simulation on large heterogeneous computing platforms where performance issues related to communication latency, load balancing and fault-tolerance are particularly exacerbated. Indeed, so far, the asynchronous iterative model constitutes the only parallel computational model which does not require global serial data management, hence, does not suffer from the Amdahl's speedup asymptotic limit~\cite{Amdahl1967}. On practical aspects, then, since global synchronization never needs to be performed, such a model provides us with a quite straightforward fault-resilience ability and less sensitivity to both communication delays and unbalanced load.

Asynchronous iterations arose as a generalization of classical fixed-point iterations, both in linear and nonlinear frameworks (see, e.g., \cite{ChazMir1969, Baudet1978}). Therefore, regarding domain decomposition methods, their application mainly targeted the overlapping Schwarz framework (see, e.g., \cite{BahiEtAl1996, FromEtAl1997, SpitEtAl2001, SpitEtAl2003}). Recently, asynchronous nonoverlapping decomposition has been addressed within the frameworks of optimized Schwarz methods (see, e.g., \cite{MagEtAl2017}), the Parareal time-decomposed time integration method (see \cite{MagGBen2018}) and a Jacobi-based sub-structuring approach (see \cite{MagVen2018}). In this paper, we address the design of an asynchronous iterative model within the primal Schur domain decomposition framework which originally features no fixed-point iterative scheme.

Precisely, considering $p$ joined nonoverlapping subdomains, we address a Schur complement inversion problem of the form
\[
\sum_{i=1}^{p} {R^{(i)}}^{\mathsf T} S^{(i)} R^{(i)} z = \sum_{i=1}^{p} {R^{(i)}}^{\mathsf T} d^{(i)},
\]
defined on the joint interface between the subdomains. To design an asynchronous fixed-point iterative scheme within such a framework, two main issues are to be handled.
First, a classical asynchronous iterative model would require a decomposition of the form
\[
z =
\begin{bmatrix}
z_{1} & \cdots & z_{p}
\end{bmatrix}^{\mathsf T},
\]
while the primal nonoverlapping domain decomposition rather induces
\[
z = z^{(1)} = \cdots = z^{(p)}.
\]
Second, it can happen that the entries of the matrices $S^{(i)}$ are not explicitly known, which would prevent us from being able to consider classical matrix splittings requiring, for instance, the diagonal entries of
\[
\sum_{i=1}^{p} {R^{(i)}}^{\mathsf T} S^{(i)} R^{(i)}.
\]
We propose here both a suitable asynchronous iterative model and practical applicable matrix splittings. General convergence conditions are provided, which consist of spectral radius bounds usually associated to asynchronous fixed-point methods. Convergent splittings of symmetric positive definite (SPD) matrices have been studied at a more general scope, e.g., in \cite{AxelKolo1994}, where it is shown that constructing a convergent splitting of an SPD matrix can be achieved by constructing a convergent splitting of a related, easily deduced, $\mathsf{M}$-matrix. A generalization of this result is due to \cite{YunKim2004}, where the deduced matrix is an $\mathsf{H}$-matrix. In \cite{CastelEtAl1998}, some splitting conditions ensure both convergence and asynchronous convergence of two-stage multi-splitting methods for Hermitian PD matrices. The approach in \cite{YunKim2004} consists of scaling up the diagonal entries of the given SPD matrix $A$ until a strictly diagonally dominant matrix $\widehat A$ is obtained, which, therefore, is an $\mathsf{H}$-matrix. Then, a convergent splitting $\widehat A = \widehat M - \widehat N$ implies a convergent splitting $A = \widehat M - N$. We therefore provide here examples of asynchronously convergent splittings of $\mathsf{H}$-matrices, and resort to diagonal scaling in numerical experiments when necessary.

The paper is organized as follows. Section \ref{sec:ai} recalls general notions involved in the convergence analysis of the new method, general asynchronous iterative models and the primal Schur domain decomposition framework. Our main results are in Section \ref{sec:aschur}, where we formulate the relaxation-based iterative scheme and establish its asynchronous convergence conditions, including a practical class of matrix splitting. Experimental results are in Section \ref{sec:experiments}, which target Poisson's and linear elasticity problems on a portion of 3D helicoid domain. The conclusions follow in Section \ref{sec:conclusions}.

\section{Preliminary notions}
\label{sec:ai}

\subsection{Convergence analysis tools}

\begin{notation}
$\mathcal A = (a_{i,j})_{i,j \in \mathbb{N}}$ denotes a finite matrix $\mathcal A$ which entry on its $i$-th row and $j$-th column is denoted as $a_{i,j}$, $i$ ranging from $1$ to its number of rows and $j$ ranging from $1$ to its number of columns. The analogous notation $x = (x_{i})_{i \in \mathbb{N}}$ is used for column vectors.
\end{notation}
\begin{definition}
\label{def:wmn}
Let $\mathcal A = (a_{i,j})_{i,j \in \mathbb{N}}$ be a square matrix and let $x = (x_{i})_{i \in \mathbb{N}}$ and $w = (w_{i})_{i \in \mathbb{N}} > 0$ be vectors with as many entries as the number of rows in $\mathcal A$. The weighted maximum norm $\|.\|_{\infty}^{w}$ is defined as
\[
\|x\|_{\infty}^{w} := \max_{i} \frac{1}{w_{i}} |x_{i}|,
\qquad
\|\mathcal A\|_{\infty}^{w} := \max_{i} \frac{1}{w_{i}} \sum_{j} |a_{i,j}| w_{j}.
\]
\end{definition}
\begin{notation}
$\rho(\mathcal A)$ denotes the spectral radius of a matrix $\mathcal A$.
\end{notation}
\begin{notation}
$|\mathcal A|$ denotes the entry-wise absolute value of a matrix $\mathcal A$.
\end{notation}
\begin{lemma}
\label{lem:pfc}
Let $\mathcal A$ be a square matrix. Then,
\[
\rho(\left|\mathcal A\right|) < 1 \quad \iff \quad \exists \ w > 0 : \left\|\mathcal A\right\|_{\infty}^{w} < 1.
\]
\end{lemma}
\proof
See, e.g., Corollary 6.1 in \cite{BertTsit1989}. % (p. 150)
\endproof
\begin{notation}
For any two matrices $\mathcal A$ and $\mathcal B$ with same numbers of rows and columns, comparisons of the type $\mathcal A < \mathcal B$ denote entry-wise comparisons and $\mathcal A \ge 0$ therefore denotes a non-negative (entry-wise) matrix.
\end{notation}
\begin{lemma}
\label{lem:varga2000_theo2.21}
For any two square matrices $\mathcal A$ and $\mathcal B$,
\[
|\mathcal A| \le \mathcal B \quad \implies \quad \rho(\mathcal A) \le \rho(\mathcal B).
\]
\end{lemma}
\proof
See, e.g., Theorem 2.21 in \cite{Varga2000}.
\endproof
\begin{definition}[$\mathsf{M}$-matrix]
\label{def:mmatrix}
An $\mathsf{M}$-matrix $\mathcal A$ is a square matrix for which there exists a real $\alpha$ such that
\[
\alpha I - \mathcal A \ge 0, \qquad \alpha > \rho(\alpha I - \mathcal A),
\]
where $I$ denotes the identity matrix.
\end{definition}
\begin{remark}[see, e.g., \cite{Fan1960}]
\label{rem:mmatrix}
An $\mathsf{M}$-matrix $\mathcal A = (a_{i,j})_{i,j \in \mathbb{N}}$ can equivalently be defined as a non-singular square matrix such that
\[
\left\{
\begin{array}{lcll}
a_{i,j} & \le & 0, & i \ne j,\\
\mathcal A^{-1} & \ge & 0.
\end{array}
\right.
\]
\end{remark}
\begin{definition}[Comparison matrix]
\label{def:compmatrix}
Let $\mathcal A = (a_{i,j})_{i,j \in \mathbb{N}}$ be a square matrix. Its comparison matrix $\langle \mathcal A \rangle = (\widetilde a_{i,j})_{i,j \in \mathbb{N}}$  is given by
\[
\left\{
\begin{array}{lcll}
\widetilde a_{i,i} & := & |a_{i,i}|, &\\
\widetilde a_{i,j} & := & -|a_{i,j}|, & i \ne j.
\end{array}
\right.
\]
\end{definition}
\begin{definition}[$\mathsf{H}$-matrix]
\label{def:hmatrix}
An $\mathsf{H}$-matrix $\mathcal A$ is a matrix such that $\langle \mathcal A \rangle$ is an $\mathsf{M}$-matrix.
\end{definition}
\begin{remark}[see, e.g.,~\cite{YunKim2004}]
According to, e.g., Theorem 5' in \cite{Fan1958}, an $\mathsf{H}$-matrix $\mathcal A = (a_{i,j})_{i,j \in \mathbb{N}}$ can equivalently be defined as a generalized strictly diagonally dominant matrix, i.e.,
\[
\exists\ u = (u_{j})_{j \in \mathbb{N}} > 0 : \quad \forall i, \ |a_{i,i}| u_{i} > \sum_{j \ne i} |a_{i,j}| u_{j}.
\]
\end{remark}
\begin{definition}[$\mathsf H$-splitting]
\label{def:hsplit}
Let $\mathcal M$ be a non-singular matrix. An $\mathsf H$-splitting $\mathcal A = \mathcal M - \mathcal N$ of a matrix $\mathcal A$ is a splitting such that $\langle \mathcal M \rangle - |\mathcal N|$ is an $\mathsf M$-matrix.
\end{definition}
\begin{lemma}
\label{lem:fs1992}
If $\mathcal A = \mathcal M - \mathcal N$ is an $\mathsf H$-splitting, then
\[
\rho(|I - \mathcal M^{-1} \mathcal A|) < 1
\]
with $I$ denoting the identity matrix.
\end{lemma}
\proof
This is straightforward from Proof of Theorem 3.4 (c) in \cite{FromSzyld1992}.
\endproof
\begin{notation}
\label{not:vmax}
Let $x = (x_{i})_{i \in \mathbb{N}}$ and $y = (y_{i})_{i \in \mathbb{N}}$ be two vectors of the same size. Then, $\max (x, y)$ denotes their entry-wise maximum, i.e., the vector with entries given by
\[
(\max (x, y))_{i} := \max \{x_{i}, y_{i}\}.
\]
\end{notation}
\begin{definition}[$|.|$-contraction]
\label{def:abs_contraction}
Let $E$ be a vector space. A $|.|$-contracting function (contracting with respect to the absolute value operator $|.|$)
\[
\mathcal F : E^{m} \to E, \qquad m \in \mathbb{N},
\]
is a function for which there exists a matrix $\mathcal T \ge 0$ with $\rho(\mathcal T) < 1$ such that
\[
\forall X, Y \in E^{m}, \qquad |\mathcal F(X) - \mathcal F(Y)| \le \mathcal T \max \left(|x^{(1)} - y^{(1)}|, \ldots, |x^{(m)} - y^{(m)}|\right)
\]
with $X := (x^{(1)}, \ldots, x^{(m)})$ and $Y := (y^{(1)}, \ldots, y^{(m)})$.
\end{definition}

\subsection{Asynchronous iterations framework}

Let
\begin{equation}
\label{eq:axeb}
A x = b
\end{equation}
be a linear system of $n$ equations, let $A = M - N$ be a splitting of $A$, and let $f$ be the associated function defined by
\begin{equation}
\label{eq:f}
f(x) = (I - M^{-1} A) x + M^{-1} b,
\end{equation}
where $I$ denotes the identity matrix. Then, one verifies
\[
A x = b \quad \iff \quad x = f(x),
\]
which allows for a fixed-point formulation of the problem \eqref{eq:axeb}. Classical fixed-point iterations generate a sequence $\{x^{k}\}_{k \in \mathbb{N}}$ such that
\begin{equation}
\label{eq:si}
x^{k+1} = f(x^{k}),
\end{equation}
which is expected to converge toward the solution $x^{*}$ of \eqref{eq:axeb}, for any given $x^{0}$. Equivalently, considering a decomposition
\[
x = 
\begin{bmatrix}
x_{1} & \cdots & x_{p}
\end{bmatrix}^{\mathsf T},
\qquad
f(x) = 
\begin{bmatrix}
f_{1}(x) & \cdots & f_{p}(x)
\end{bmatrix}^{\mathsf T},
\qquad
p \le n,
\]
the iterative model \eqref{eq:si} can be reformulated as
\[
x_{i}^{k+1} = f_{i}(x_{1}^{k}, \ldots, x_{p}^{k}) \qquad \forall i \in \{1, \ldots, p\}.
\]
Asynchronous iterations arose from the free steering model which allows for considering that a global iteration consists of any update of an arbitrary subset
\[
P^{(k)} \subseteq \{1, \ldots, p\}
\]
of components of $x^{k}$ (see, e.g., \cite{Schechter1959}). Additionally, if one no longer explicitly performs such global iterations, then any component may be updated concurrently to its read access, and therefore, any update $x_{i}^{k+1}$ is potentially based on outdated components $x_{j}^{\tau_{j}^{i}(k)}$, $j \in \{1, \ldots, p\}$, i.e., with
\[
\tau_{j}^{i}(k) \le k.
\]
The corresponding asynchronous iterative model is thus given by (see, e.g., \cite{BertTsit1991})
\begin{equation}
\label{eq:ai}
x_{i}^{k+1} = \left \{
\begin{array}{ll}
f_{i}\left(x_{1}^{\tau_{1}^{i}(k)}, \ldots, x_{p}^{\tau_{p}^{i}(k)}\right) & \forall i \in P^{(k)},\\
x_{i}^{k} & \forall i \notin P^{(k)},
\end{array}
\right.
\end{equation}
where it is however assumed that every component is infinitely often updated, i.e.,
\begin{equation}
\label{eq:ai_ass1}
\forall i \in \{1, \dots, p\}, \quad \operatorname{card}\{k \in \mathbb{N} \ | \ i \in P^{(k)}\} = +\infty,
\end{equation}
and that access delays to updated components are bounded, so as to have
\begin{equation}
\label{eq:ai_ass2}
\forall i, j \in \{1, \dots, p\}, \quad \lim_{k \to +\infty} \tau_{j}^{i}(k) = +\infty.
\end{equation}
In \cite{Baudet1978}, nonlinear fixed-point problems of the form
\[
x = \widetilde f(x, x, \ldots, x), \qquad\qquad \widetilde f : E^{m} \to E, \quad m \in \mathbb{N},
\]
are introduced to analyze Newton asynchronous iterations. It yielded a generalization of \eqref{eq:ai} into
\begin{equation}
\label{eq:gai}
x_{i}^{k+1} = \left \{
\begin{array}{ll}
\widetilde f_{i}\left(x_{1}^{\tau_{1,1}^{i}(k)}, \ldots, x_{p}^{\tau_{p,1}^{i}(k)}, \ldots, x_{1}^{\tau_{1,m}^{i}(k)}, \ldots, x_{p}^{\tau_{p,m}^{i}(k)}\right) & \forall i \in P^{(k)},\\
x_{i}^{k} & \forall i \notin P^{(k)},
\end{array}
\right.
\end{equation}
so that \eqref{eq:ai} corresponds to the case $m = 1$.
\begin{theorem}[Chazan and Miranker (1969)]
\label{theo:cm1969}
The asynchronous iterative model \eqref{eq:ai} is convergent for any $x^{0}$, $\{P^{(k)}\}_{k \in \mathbb{N}}$ and $\tau_{j}^{i}$ with $i, j \in \{1, \ldots, p\}$ if and only if
\[
\rho(|I - M^{-1} A|) < 1.
\]
\end{theorem}
\proof
See \cite{ChazMir1969}.
\endproof
\begin{theorem}[Baudet (1978)]
\label{theo:b1978}
The asynchronous iterative model \eqref{eq:gai} is convergent if $\widetilde f$ is $|.|$-contracting.
\end{theorem}
\proof
See \cite{Baudet1978}.
\endproof

\subsection{Domain decomposition framework}

Let
\begin{equation}
\label{eq:ddm:axeb}
\begin{bmatrix}
A_{II}^{(1)} & 0 & A_{I\Gamma}^{(1)}\\
0 & A_{II}^{(2)} & A_{I\Gamma}^{(2)}\\
A_{\Gamma I}^{(1)} & A_{\Gamma I}^{(2)} & A_{\Gamma\Gamma}^{(1)}+A_{\Gamma\Gamma}^{(2)}
\end{bmatrix}
\begin{bmatrix}
x_{I}^{(1)}\\
x_{I}^{(2)}\\
x_{\Gamma}^{ }
\end{bmatrix}
=
\begin{bmatrix}
b_{I}^{(1)}\\
b_{I}^{(2)}\\
b_{\Gamma}^{(1)}+b_{\Gamma}^{(2)}
\end{bmatrix}
\end{equation}
be a linear problem defined on two nonoverlapping (sub)domains $\Omega_{1}$ and $\Omega_{2}$ joined by an interface $\Gamma$ with respective solutions $x_{I}^{(1)}$, $x_{I}^{(2)}$ and $x_{\Gamma}$. By eliminating either $x_{I}^{(2)}$ or $x_{I}^{(1)}$, it yields either
\[
\begin{bmatrix}
A_{II}^{(1)} & A_{I\Gamma}^{(1)}\\
A_{\Gamma I}^{(1)} & A_{\Gamma\Gamma}^{(1)} + A_{\Gamma\Gamma}^{(2)} - A_{\Gamma I}^{(2)} {A_{II}^{(2)}}^{-1} A_{I\Gamma}^{(2)}
\end{bmatrix}
\begin{bmatrix}
x_{I}^{(1)}\\
x_{\Gamma}^{ }
\end{bmatrix}
=
\begin{bmatrix}
b_{I}^{(1)}\\
b_{\Gamma}^{(1)}+b_{\Gamma}^{(2)} - A_{\Gamma I}^{(2)} {A_{II}^{(2)}}^{-1} b_{I}^{(2)}
\end{bmatrix}
\]
or
\[
\begin{bmatrix}
A_{II}^{(2)} & A_{I\Gamma}^{(2)}\\
A_{\Gamma I}^{(2)} & A_{\Gamma\Gamma}^{(1)} + A_{\Gamma\Gamma}^{(2)} - A_{\Gamma I}^{(1)} {A_{II}^{(1)}}^{-1} A_{I\Gamma}^{(1)}
\end{bmatrix}
\begin{bmatrix}
x_{I}^{(2)}\\
x_{\Gamma}^{ }
\end{bmatrix}
=
\begin{bmatrix}
b_{I}^{(2)}\\
b_{\Gamma}^{(1)}+b_{\Gamma}^{(2)} - A_{\Gamma I}^{(1)} {A_{II}^{(1)}}^{-1} b_{I}^{(1)}
\end{bmatrix},
\]
which are two ways of solving the same initial problem. Then, from each of these two parallel points of view $i \in \{1,2\}$, and considering that data related to $j \in \{1,2\}$, $j \ne i$, are unknown, one sets up two different equations
\[
\begin{bmatrix}
A_{II}^{(1)} & A_{I\Gamma}^{(1)}\\
A_{\Gamma I}^{(1)} & A_{\Gamma\Gamma}^{(1)} + \Lambda_{\Gamma\Gamma}^{(1)}
\end{bmatrix}
\begin{bmatrix}
x_{I}^{(1)}\\
x_{\Gamma}^{(1)}
\end{bmatrix}
=
\begin{bmatrix}
b_{I}^{(1)}\\
b_{\Gamma}^{(1)}+\lambda_{\Gamma}^{(1)}
\end{bmatrix}
\]
and
\[
\begin{bmatrix}
A_{II}^{(2)} & A_{I\Gamma}^{(2)}\\
A_{\Gamma I}^{(2)} & A_{\Gamma\Gamma}^{(2)} + \Lambda_{\Gamma\Gamma}^{(2)}
\end{bmatrix}
\begin{bmatrix}
x_{I}^{(2)}\\
x_{\Gamma}^{(2)}
\end{bmatrix}
=
\begin{bmatrix}
b_{I}^{(2)}\\
b_{\Gamma}^{(2)} + \lambda_{\Gamma}^{(2)}
\end{bmatrix},
\]
which solutions, for any choice of $\Lambda_{\Gamma\Gamma}^{(1)}$ and $\Lambda_{\Gamma\Gamma}^{(2)}$, are those of the global problem \eqref{eq:ddm:axeb} if and only if (see, e.g., Theorem~2.1 in \cite{MagoulesEtAl2004})
\[
\left\{
\begin{aligned}
x_{\Gamma}^{(1)} & = x_{\Gamma}^{(2)},\\
\lambda_{\Gamma}^{(1)} - \Lambda_{\Gamma\Gamma}^{(1)} x_{\Gamma}^{(1)} & = -\left(\lambda_{\Gamma}^{(2)} - \Lambda_{\Gamma\Gamma}^{(2)} x_{\Gamma}^{(2)}\right).
\end{aligned}
\right.
\]
Then, again, by eliminating $x_{I}^{(1)}$ and $x_{I}^{(2)}$, one reaches an equivalent problem defined on $\Gamma$ by
\[
\left\{
\begin{aligned}
\left(A_{\Gamma\Gamma}^{(1)} - A_{\Gamma I}^{(1)} {A_{II}^{(1)}}^{-1} A_{I\Gamma}^{(1)} + \Lambda_{\Gamma\Gamma}^{(1)}\right) x_{\Gamma}^{(1)} & = b_{\Gamma}^{(1)} + \lambda_{\Gamma}^{(1)} - A_{\Gamma I}^{(1)} {A_{II}^{(1)}}^{-1} b_{I}^{(1)},\\
\left(A_{\Gamma\Gamma}^{(2)} - A_{\Gamma I}^{(2)} {A_{II}^{(2)}}^{-1} A_{I\Gamma}^{(2)} + \Lambda_{\Gamma\Gamma}^{(2)}\right) x_{\Gamma}^{(2)} & = b_{\Gamma}^{(2)} + \lambda_{\Gamma}^{(2)} - A_{\Gamma I}^{(2)} {A_{II}^{(2)}}^{-1} b_{I}^{(2)},\\
x_{\Gamma}^{(1)} & = x_{\Gamma}^{(2)},\\
\lambda_{\Gamma}^{(1)} - \Lambda_{\Gamma\Gamma}^{(1)} x_{\Gamma}^{(1)} & = -\left(\lambda_{\Gamma}^{(2)} - \Lambda_{\Gamma\Gamma}^{(2)} x_{\Gamma}^{(2)}\right).
\end{aligned}
\right.
\]
Let us finally identify the local Schur complements
\[
S_{\Gamma\Gamma}^{(i)} := A_{\Gamma\Gamma}^{(i)} - A_{\Gamma I}^{(i)} {A_{II}^{(i)}}^{-1} A_{I\Gamma}^{(i)}, \qquad i \in \{1,2\},
\]
and define
\[
d_{\Gamma}^{(i)} := b_{\Gamma}^{(i)} - A_{\Gamma I}^{(i)} {A_{II}^{(i)}}^{-1} b_{I}^{(i)}, \qquad i \in \{1,2\}.
\]
Then, eliminating the unknowns $\lambda_{\Gamma}^{(i)}$, $i \in \{1,2\}$, leads to
\[
\left\{
\begin{aligned}
x_{\Gamma}^{(1)} & = x_{\Gamma}^{(2)},\\
S_{\Gamma\Gamma}^{(1)} x_{\Gamma}^{(1)} + S_{\Gamma\Gamma}^{(2)} x_{\Gamma}^{(2)} & = d_{\Gamma}^{(1)} + d_{\Gamma}^{(2)},
\end{aligned}
\right.
\]
and, hence, to each of the identical equations
\begin{equation}
\label{eq:sxed}
\left(S_{\Gamma\Gamma}^{(1)} + S_{\Gamma\Gamma}^{(2)}\right) x_{\Gamma}^{(1)} = d_{\Gamma}^{(1)} + d_{\Gamma}^{(2)},
\qquad
\left(S_{\Gamma\Gamma}^{(1)} + S_{\Gamma\Gamma}^{(2)}\right) x_{\Gamma}^{(2)} = d_{\Gamma}^{(1)} + d_{\Gamma}^{(2)},
\end{equation}
where
\[
S_{\Gamma\Gamma} := \left(S_{\Gamma\Gamma}^{(1)} + S_{\Gamma\Gamma}^{(2)}\right)
\]
is the global Schur complement from \eqref{eq:ddm:axeb}.

To generalize \eqref{eq:ddm:axeb} to $p$ subdomains, Boolean matrices $R_{\Gamma_{l}\Gamma}^{(i)}$, $i \in \{1, \ldots, p\}$, are considered to denote mappings between a whole interface $\Gamma$ and restricted parts $\Gamma_{l}$ that subdomains are locally acting on. We would then actually have
\[
A_{I\Gamma}^{(i)} = A_{I\Gamma_{l}}^{(i)} R_{\Gamma_{l}\Gamma}^{(i)},
\quad
A_{\Gamma I}^{(i)} = R_{\Gamma\Gamma_{l}}^{(i)} A_{\Gamma_{l} I}^{(i)},
\quad
A_{\Gamma\Gamma}^{(i)} = R_{\Gamma\Gamma_{l}}^{(i)} A_{\Gamma_{l}\Gamma_{l}}^{(i)} R_{\Gamma_{l}\Gamma}^{(i)},
\quad
b_{\Gamma}^{(i)} = R_{\Gamma\Gamma_{l}}^{(i)} b_{\Gamma_{l}}^{(i)}
\]
with
\[
R_{\Gamma\Gamma_{l}}^{(i)} = {R_{\Gamma_{l}\Gamma}^{(i)}}^{\mathsf T},
\]
and, further,
\[
S_{\Gamma\Gamma}^{(i)} = R_{\Gamma\Gamma_{l}}^{(i)} S_{\Gamma_{l}\Gamma_{l}}^{(i)} R_{\Gamma_{l}\Gamma}^{(i)},
\qquad
d_{\Gamma}^{(i)} = R_{\Gamma\Gamma_{l}}^{(i)} d_{\Gamma_{l}}^{(i)}
\]
with
\[
S_{\Gamma_{l}\Gamma_{l}}^{(i)} := A_{\Gamma_{l}\Gamma_{l}}^{(i)} - A_{\Gamma_{l} I}^{(i)} {A_{II}^{(i)}}^{-1} A_{I\Gamma_{l}}^{(i)},
\qquad
d_{\Gamma_{l}}^{(i)} := b_{\Gamma_{l}}^{(i)} - A_{\Gamma_{l} I}^{(i)} {A_{II}^{(i)}}^{-1} b_{I}^{(i)}.
\]
The interface problem \eqref{eq:sxed} is therefore generalized as
\[
S_{\Gamma\Gamma}^{} x_{\Gamma}^{} = d_{\Gamma}^{}
\]
with
\[
S_{\Gamma\Gamma}^{} := \sum_{i=1}^{p} S_{\Gamma\Gamma}^{(i)}, \qquad d_{\Gamma}^{} := \sum_{i=1}^{p} d_{\Gamma}^{(i)}.
\]

\section{Asynchronous primal Schur method}
\label{sec:aschur}

\subsection{Iterative scheme}

Let us consider, as in \eqref{eq:sxed}, $p$ identical equations
\[
\sum_{j=1}^{p} S_{\Gamma\Gamma}^{(j)} x_{\Gamma}^{(i)} = \sum_{j=1}^{p} d_{\Gamma}^{(j)}, \qquad i \in \{1, \ldots, p\}.
\]
A splitting $S_{\Gamma\Gamma} = M_{\Gamma\Gamma} - N_{\Gamma\Gamma}$ of each of them leads to $p$ identical iterative procedures
\[
x_{\Gamma}^{(i),k+1} = \left(I_{\Gamma\Gamma} - M_{\Gamma\Gamma}^{-1} \sum_{j=1}^{p} S_{\Gamma\Gamma}^{(j)}\right) x_{\Gamma}^{(i),k} + M_{\Gamma\Gamma}^{-1} \sum_{j=1}^{p} d_{\Gamma}^{(j)}, \qquad i \in \{1, \ldots, p\}
\]
with $I_{\Gamma\Gamma}$ denoting the identity matrix defined on $\Gamma$. Let us, then, consider $p$ matrices $I_{\Gamma\Gamma}^{(i)}$, $i \in \{1, \ldots, p\}$, such that
\[
\sum_{i=1}^{p} I_{\Gamma\Gamma}^{(i)} = I_{\Gamma\Gamma}.
\]
We thus have
\[
x_{\Gamma}^{(i),k+1} = \left(\sum_{j=1}^{p} I_{\Gamma\Gamma}^{(j)} - M_{\Gamma\Gamma}^{-1} S_{\Gamma\Gamma}^{(j)}\right) x_{\Gamma}^{(i),k} + M_{\Gamma\Gamma}^{-1} \sum_{j=1}^{p} d_{\Gamma}^{(j)}, \qquad i \in \{1, \ldots, p\}.
\]
Then, as these are equivalent iterations, we also have
\[
x_{\Gamma}^{(1),k} = \cdots = x_{\Gamma}^{(p),k} \qquad \forall k,
\]
which allows us to replace $x_{\Gamma}^{(i),k}$ by any $x_{\Gamma}^{(j),k}$ and get
\[
\begin{aligned}
x_{\Gamma}^{(i),k+1} & = \sum_{j=1}^{p} \left(I_{\Gamma\Gamma}^{(j)} - M_{\Gamma\Gamma}^{-1} S_{\Gamma\Gamma}^{(j)}\right) x_{\Gamma}^{(j),k} + M_{\Gamma\Gamma}^{-1} \sum_{j=1}^{p} d_{\Gamma}^{(j)}\\
& = \sum_{j=1}^{p} I_{\Gamma\Gamma}^{(j)} x_{\Gamma}^{(j),k} + M_{\Gamma\Gamma}^{-1} \left(d_{\Gamma}^{(j)} - S_{\Gamma\Gamma}^{(j)} x_{\Gamma}^{(j),k}\right).
\end{aligned}
\]
The corresponding parallel iterative scheme is thus given by
\begin{equation}
\label{eq:ddm:si}
\left\{
\begin{array}{lcl}
y_{\Gamma}^{(i),k} & := & I_{\Gamma\Gamma}^{(i)} x_{\Gamma}^{(i),k} + M_{\Gamma\Gamma}^{-1} \left(d_{\Gamma}^{(i)} - S_{\Gamma\Gamma}^{(i)} x_{\Gamma}^{(i),k}\right),\\
x_{\Gamma}^{(i),k+1} & = & \displaystyle\sum_{j=1}^{p} y_{\Gamma}^{(j),k} \quad \forall i \in \{1, \ldots, p\},
\end{array}
\right.
\end{equation}
which exhibits communication only for the vectors $y_{\Gamma}^{(j),k}$ with $j \ne i$. This therefore generalizes to an asynchronous iterative scheme
\begin{equation}
\label{eq:ddm:ai}
\left\{
\begin{array}{lcl}
y_{\Gamma}^{(i),k} & := & I_{\Gamma\Gamma}^{(i)} x_{\Gamma}^{(i),k} + M_{\Gamma\Gamma}^{-1} \left(d_{\Gamma}^{(i)} - S_{\Gamma\Gamma}^{(i)} x_{\Gamma}^{(i),k}\right),\\
x_{\Gamma}^{(i),k+1} & = & \left\{
\begin{array}{ll}
\displaystyle\sum_{j=1}^{p} y_{\Gamma}^{(j),\tau_{j}^{i}(k)} & \forall i \in P^{(k)},\\
x_{\Gamma}^{(i),k} & \forall i \notin P^{(k)},
\end{array}
\right.
\end{array}
\right.
\end{equation}
which lies in the framework of the general model \eqref{eq:gai} with, for each process $i \in \{1, \ldots, p\}$, $\widetilde f := f^{(i)}$ being given by
\begin{equation}
\label{eq:ddm:f}
f^{(i)}(x_{\Gamma}^{(1)}, \ldots, x_{\Gamma}^{(p)}) = \sum_{j=1}^{p} I_{\Gamma\Gamma}^{(j)} x_{\Gamma}^{(j)} + M_{\Gamma\Gamma}^{-1} \left(d_{\Gamma}^{(j)} - S_{\Gamma\Gamma}^{(j)} x_{\Gamma}^{(j)}\right),
\end{equation}
and verifying
\[
S_{\Gamma\Gamma}^{} x_{\Gamma}^{} = d_{\Gamma}^{}
\quad \iff \quad
x_{\Gamma}^{ } = f^{(i)}(x_{\Gamma}^{ }, \ldots, x_{\Gamma}^{ }).
\]

While the first equation in both \eqref{eq:ddm:si} and \eqref{eq:ddm:ai} is defined on the whole interface $\Gamma$, in practice, it can be operated in each local subspace $\Gamma_{l}$ if, for instance, $M_{\Gamma\Gamma}$ is diagonal. Indeed, if $M_{\Gamma\Gamma}$ is diagonal, then one satisfies
\[
\begin{aligned}
y_{\Gamma}^{(i),k} & = I_{\Gamma\Gamma}^{(i)} x_{\Gamma}^{(i),k} + M_{\Gamma\Gamma}^{-1} \left(d_{\Gamma}^{(i)} - S_{\Gamma\Gamma}^{(i)} x_{\Gamma}^{(i),k}\right)\\
& = R_{\Gamma\Gamma_{l}}^{(i)} I_{\Gamma_{l}\Gamma_{l}}^{(i)} R_{\Gamma_{l}\Gamma}^{(i)} x_{\Gamma}^{(i),k} + M_{\Gamma\Gamma}^{-1} \left(R_{\Gamma\Gamma_{l}}^{(i)} d_{\Gamma_{l}}^{(i)} - R_{\Gamma\Gamma_{l}}^{(i)} S_{\Gamma_{l}\Gamma_{l}}^{(i)} R_{\Gamma_{l}\Gamma}^{(i)} x_{\Gamma}^{(i),k}\right)\\
& = R_{\Gamma\Gamma_{l}}^{(i)} I_{\Gamma_{l}\Gamma_{l}}^{(i)} x_{\Gamma_{l}}^{(i),k} + M_{\Gamma\Gamma}^{-1} \left(R_{\Gamma\Gamma_{l}}^{(i)} d_{\Gamma_{l}}^{(i)} - R_{\Gamma\Gamma_{l}}^{(i)} S_{\Gamma_{l}\Gamma_{l}}^{(i)} x_{\Gamma_{l}}^{(i),k}\right)\\
& = R_{\Gamma\Gamma_{l}}^{(i)} I_{\Gamma_{l}\Gamma_{l}}^{(i)} x_{\Gamma_{l}}^{(i),k} + M_{\Gamma\Gamma}^{-1} R_{\Gamma\Gamma_{l}}^{(i)} \left(d_{\Gamma_{l}}^{(i)} - S_{\Gamma_{l}\Gamma_{l}}^{(i)} x_{\Gamma_{l}}^{(i),k}\right)\\
& = R_{\Gamma\Gamma_{l}}^{(i)} I_{\Gamma_{l}\Gamma_{l}}^{(i)} x_{\Gamma_{l}}^{(i),k} + R_{\Gamma\Gamma_{l}}^{(i)} R_{\Gamma_{l}\Gamma}^{(i)} M_{\Gamma\Gamma}^{-1} R_{\Gamma\Gamma_{l}}^{(i)} \left(d_{\Gamma_{l}}^{(i)} - S_{\Gamma_{l}\Gamma_{l}}^{(i)} x_{\Gamma_{l}}^{(i),k}\right)\\
& = R_{\Gamma\Gamma_{l}}^{(i)} \left[I_{\Gamma_{l}\Gamma_{l}}^{(i)} x_{\Gamma_{l}}^{(i),k} + R_{\Gamma_{l}\Gamma}^{(i)} M_{\Gamma\Gamma}^{-1} R_{\Gamma\Gamma_{l}}^{(i)} \left(d_{\Gamma_{l}}^{(i)} - S_{\Gamma_{l}\Gamma_{l}}^{(i)} x_{\Gamma_{l}}^{(i),k}\right)\right],
\end{aligned}
\]
which leads to
\[
\left\{
\begin{array}{lcl}
y_{\Gamma_{l}}^{(i),k} & := & I_{\Gamma_{l}\Gamma_{l}}^{(i)} x_{\Gamma_{l}}^{(i),k} + R_{\Gamma_{l}\Gamma}^{(i)} M_{\Gamma\Gamma}^{-1} R_{\Gamma\Gamma_{l}}^{(i)} \left(d_{\Gamma_{l}}^{(i)} - S_{\Gamma_{l}\Gamma_{l}}^{(i)} x_{\Gamma_{l}}^{(i),k}\right),\\
x_{\Gamma}^{(i),k+1} & = & \displaystyle\sum_{j=1}^{p} R_{\Gamma\Gamma_{l}}^{(j)} y_{\Gamma_{l}}^{(j),k}.
\end{array}
\right.
\]

\subsection{Convergence conditions}

\begin{lemma}
\label{lem:abs_contraction}
Let
\[
\mathcal F : E^{m \times n} \mapsto E^{n},
\qquad
\mathcal F :=
\begin{bmatrix}
\mathcal F^{(1)} & \cdots & \mathcal F^{(n)}
\end{bmatrix}^{\mathsf T},
\qquad
n, m \in \mathbb{N},
\]
be a function derived from $n$ given functions
\[
\mathcal F^{(i)} : E^{m} \mapsto E, \qquad i \in \{1, \ldots, n\}.
\]
Then, $\mathcal F$ is $|.|$-contracting if $\mathcal F^{(i)}$ is $|.|$-contracting for all $i \in \{1, \ldots, n\}$.
\end{lemma}
\proof
We conduct the proof for $n = 2$ and $m = 2$. The generalization to any $n$ and $m$ follows in the same way.

According to Definition \ref{def:abs_contraction}, let $\mathcal T^{(1)}$ and $\mathcal T^{(2)}$ be the matrices of contraction of $\mathcal F^{(1)}$ and $\mathcal F^{(2)}$ respectively, which satisfy
\[
\mathcal T^{(1)} \ge 0, \quad \rho(\mathcal T^{(1)}) < 1, \qquad \mathcal T^{(2)} \ge 0, \quad \rho(\mathcal T^{(2)}) < 1.
\]
Let us consider two arbitrary vectors
\[
X := (x^{(1,1)}, x^{(1,2)}, x^{(2,1)}, x^{(2,2)}), \qquad Y := (y^{(1,1)}, y^{(1,2)}, y^{(2,1)}, y^{(2,2)})
\]
with $x^{(i,j)}, y^{(i,j)} \in E$, where $i,j \in \{1,2\}$. Then, we have
\[
\begin{aligned}
|\mathcal F(X) - \mathcal F(Y)| & =
\begin{bmatrix}
|\mathcal F^{(1)}(x^{(1,1)}, x^{(1,2)}) - \mathcal F^{(1)}(y^{(1,1)}, y^{(1,2)})|\\
|\mathcal F^{(2)}(x^{(2,1)}, x^{(2,2)}) - \mathcal F^{(2)}(y^{(2,1)}, y^{(2,2)})|
\end{bmatrix}\\
& \le
\begin{bmatrix}
\mathcal T^{(1)} \max \left(|x^{(1,1)} - y^{(1,1)}|, |x^{(1,2)} - y^{(1,2)}|\right)\\
\mathcal T^{(2)} \max \left(|x^{(2,1)} - y^{(2,1)}|, |x^{(2,2)} - y^{(2,2)}|\right)
\end{bmatrix}\\
& =
\begin{bmatrix}
\mathcal T^{(1)} & 0\\
0 & \mathcal T^{(2)}
\end{bmatrix}
\begin{bmatrix}
\max \left(|x^{(1,1)} - y^{(1,1)}|, |x^{(1,2)} - y^{(1,2)}|\right)\\
\max \left(|x^{(2,1)} - y^{(2,1)}|, |x^{(2,2)} - y^{(2,2)}|\right)
\end{bmatrix}\\
& =
\begin{bmatrix}
\mathcal T^{(1)} & 0\\
0 & \mathcal T^{(2)}
\end{bmatrix}
\max \left(\left|
\begin{bmatrix}
x^{(1,1)}\\
x^{(2,1)}
\end{bmatrix}
-
\begin{bmatrix}
y^{(1,1)}\\
y^{(2,1)}
\end{bmatrix}
\right|, \left|
\begin{bmatrix}
x^{(1,2)}\\
x^{(2,2)}
\end{bmatrix}
-
\begin{bmatrix}
y^{(1,2)}\\
y^{(2,2)}
\end{bmatrix}
\right|\right),
\end{aligned}
\]
and, hence,
\[
\mathcal T \ge 0, \quad \rho(\mathcal T) < 1, \qquad 
\mathcal T =
\begin{bmatrix}
\mathcal T^{(1)} & 0\\
0 & \mathcal T^{(2)}
\end{bmatrix},
\]
which thus concludes the proof.
\endproof
\begin{theorem}
\label{theo:ddm:ai}
The asynchronous primal Schur method \eqref{eq:ddm:ai} is convergent if
\begin{equation}
\label{eq:ddm:ai_cond}
\rho\left(\sum_{i=1}^{p} \left|I_{\Gamma\Gamma}^{(i)} - M_{\Gamma\Gamma}^{-1} S_{\Gamma\Gamma}^{(i)}\right|\right) < 1.
\end{equation}
\end{theorem}
\proof
Let us consider arbitrary vectors $x_{\Gamma}^{(j)}$, $y_{\Gamma}^{(j)}$ with $j \in \{1, \ldots, p\}$, defined on $\Gamma$, and
\[
X_{\Gamma^{p}} := (x_{\Gamma}^{(1)}, \ldots, x_{\Gamma}^{(p)}), \qquad Y_{\Gamma^{p}} := (y_{\Gamma}^{(1)}, \ldots, y_{\Gamma}^{(p)}).
\]
Then, according to \eqref{eq:ddm:f}, we have, for any $i \in \{1, \ldots, p\}$,
\[
\begin{aligned}
\left|f^{(i)}(X_{\Gamma^{p}}) - f^{(i)}(Y_{\Gamma^{p}})\right| & = \left|\sum_{j=1}^{p} I_{\Gamma\Gamma}^{(j)} \left(x_{\Gamma}^{(j)} - y_{\Gamma}^{(j)}\right) - M_{\Gamma\Gamma}^{-1} S_{\Gamma\Gamma}^{(j)} \left(x_{\Gamma}^{(j)} - y_{\Gamma}^{(j)}\right)\right|\\
& \le \sum_{j=1}^{p} \left|I_{\Gamma\Gamma}^{(j)} - M_{\Gamma\Gamma}^{-1} S_{\Gamma\Gamma}^{(j)}\right| \left|x_{\Gamma}^{(j)} - y_{\Gamma}^{(j)}\right|\\
& \le \sum_{j=1}^{p} \left|I_{\Gamma\Gamma}^{(j)} - M_{\Gamma\Gamma}^{-1} S_{\Gamma\Gamma}^{(j)}\right| \max_{j=1}^{p} \left|x_{\Gamma}^{(j)} - y_{\Gamma}^{(j)}\right|.
\end{aligned}
\]
Given \eqref{eq:ddm:ai_cond}, each $f^{(i)}$ is therefore $|.|$-contracting, and so, Lemma \ref{lem:abs_contraction} applies, which concludes the proof, according to Theorem \ref{theo:b1978}.
\endproof
\begin{assumption}
\label{ass:ddm:ai}
Let the matrices $I_{\Gamma\Gamma}^{(i)} - M_{\Gamma\Gamma}^{-1} S_{\Gamma\Gamma}^{(i)} = (q^{(i)}_{l,t})_{l,t \in \mathbb{N}}$, $i \in \{1, \ldots, p\}$, be entry-wise of non-different signs, i.e.,
\[
q^{(i)}_{l,t} < 0 \ \implies \ q^{(j)}_{l,t} \le 0, \qquad q^{(i)}_{l,t} > 0 \ \implies \ q^{(j)}_{l,t} \ge 0, \qquad i,j \in \{1, \ldots, p\}.
\]
\end{assumption}
\begin{corollary}
\label{cor:ddm:ai}
Under Assumption \ref{ass:ddm:ai}, the asynchronous primal Schur method \eqref{eq:ddm:ai} is convergent if
\[
\rho\left(\left|I_{\Gamma\Gamma} - M_{\Gamma\Gamma}^{-1} S_{\Gamma\Gamma}^{}\right|\right) < 1.
\]
\end{corollary}
\proof
Under Assumption \ref{ass:ddm:ai}, we have
\[
\sum_{i=1}^{p} \left|I_{\Gamma\Gamma}^{(i)} - M_{\Gamma\Gamma}^{-1} S_{\Gamma\Gamma}^{(i)}\right| = \left|\sum_{i=1}^{p} I_{\Gamma\Gamma}^{(i)} - M_{\Gamma\Gamma}^{-1} S_{\Gamma\Gamma}^{(i)}\right| = \left|I_{\Gamma\Gamma} - M_{\Gamma\Gamma}^{-1} S_{\Gamma\Gamma}^{}\right|,
\]
hence, Theorem \ref{theo:ddm:ai} directly applies.
\endproof

\subsection{Some practical convergence cases}

\begin{notation}
\label{not:cwmwmn}
Let $\mathcal A = (a_{i,j})_{i,j \in \mathbb{N}}$ be a rectangular matrix and let $w = (w_{j})_{j \in \mathbb{N}}$ and $v = (v_{i})_{i \in \mathbb{N}} > 0$ be two vectors having as many entries as, respectively, the number of columns and the number of rows in $\mathcal A$. We denote by $|\mathcal A|^{w}_{v}$ the vector with entries given by the weighted row-sums
\[
(|\mathcal A|^{w}_{v})_{i} := \frac{1}{v_{i}} \sum_{j} |a_{i,j}| w_{j}.
\]
\end{notation}
\begin{notation}
\label{not:cwvd}
Let $w = (w_{i})_{i \in \mathbb{N}}$ and $v = (v_{i})_{i \in \mathbb{N}} > 0$ be two vectors. We denote by $w/v$ the entry-wise division of $w$ by $v$, i.e., the vector with entries given by
\[
(w/v)_{i} := \frac{w_{i}}{v_{i}}.
\]
\end{notation}
\begin{remark}
\label{rem:cwmwmn}
Notation \ref{not:cwmwmn} and Notation \ref{not:cwvd} induce
\[
\|\mathcal A\|_{\infty}^{w} = \max_{i} |\mathcal A|^{w}_{w},
\qquad
|\mathcal A|^{w}_{v} = (|\mathcal A| w) / v = \sum_{j} w_{j} |\mathcal A_{j}|/v,
\]
where $\mathcal A_{j}$ denotes the $j$-th column of $\mathcal A$.
\end{remark}
\begin{lemma}
\label{lem:cwmwmn}
Let $\mathcal A = (a_{i,j})_{i,j \in \mathbb{N}}$ and $\mathcal B = (b_{j,l})_{j,l \in \mathbb{N}}$ be two matrices such that the number of columns in $\mathcal A$ equals the number of rows in $\mathcal B$. Let $z = (z_{j})_{j \in \mathbb{N}} > 0$, $v = (v_{i})_{i \in \mathbb{N}} > 0$ and $w = (w_{l})_{l \in \mathbb{N}}$ be three vectors having as many entries as, respectively, the number of columns in $\mathcal A$, the number of rows in $\mathcal A$ and the number of columns in $\mathcal B$. Let, at last, $u = (u_{j})_{j \in \mathbb{N}}$ be the vector with as many entries as the number of rows in $\mathcal B$ and given by
$
u_{j} := 1 \ \forall j.
$
Then, we have
\[
|\mathcal B|_{z}^{w} < u \quad \implies \quad |\mathcal A \mathcal B|_{v}^{w} < |\mathcal A|_{v}^{z}.
\]
\end{lemma}
\proof
According to Remark \ref{rem:cwmwmn}, we have
\[
|\mathcal A \mathcal B|_{v}^{w} = \sum_{l} w_{l} |(\mathcal A \mathcal B)_{l}|/v,
\]
where $(\mathcal A \mathcal B)_{l}$ is the $l$-th column of the product $\mathcal A \mathcal B$. Since the entries of the matrix $\mathcal A \mathcal B = (c_{i,l})_{i,l \in \mathbb{N}}$ are given by
\[
c_{i,l} := \sum_{j} a_{i,j} b_{j,l},
\]
we can also write
\[
(\mathcal A \mathcal B)_{l} = \sum_{j} b_{j,l} \mathcal A_{j},
\]
where $\mathcal A_{j}$ is the $j$-th column of $\mathcal A$. These imply that
\[
\begin{aligned}
|\mathcal A \mathcal B|_{v}^{w} & = \sum_{l} w_{l}^{} \left|\sum_{j} b_{j,l} \mathcal A_{j}\right|/v\\
& \le \sum_{l} w_{l}^{} \sum_{j} |b_{j,l} \mathcal A_{j}|/v\\
& = \sum_{l} \sum_{j} w_{l}^{} |b_{j,l}| |\mathcal A_{j}|/v\\
& = \sum_{l} \sum_{j} w_{l}^{} \frac{|b_{j,l}|}{z_{j}} z_{j} |\mathcal A_{j}|/v\\
& = \sum_{j} \left(\sum_{l} w_{l}^{} \frac{|b_{j,l}|}{z_{j}}\right) z_{j} |\mathcal A_{j}|/v.
\end{aligned}
\]
Referring to Notation \ref{not:cwmwmn}, which gives the entries of the vector $|\mathcal B|_{z}^{w}$ as
\[
(|\mathcal B|_{z}^{w})_{j} := \sum_{l} w_{l}^{} \frac{|b_{j,l}|}{z_{j}},
\]
we can further write
\[
|\mathcal A \mathcal B|_{v}^{w} \le \sum_{j} (|\mathcal B|_{z}^{w})_{j} z_{j} |\mathcal A_{j}|/v,
\]
which implies, for the entries of the vectors $|\mathcal A \mathcal B|_{v}^{w}$ and $|\mathcal A_{j}|/v$,
\[
(|\mathcal A \mathcal B|_{v}^{w})_{i} \le \sum_{j} (|\mathcal B|_{z}^{w})_{j} z_{j} (|\mathcal A_{j}|/v)_{i}.
\]
If, therefore,
$
|\mathcal B|_{z}^{w} < u,
$
then we have, still considering vectors entries,
\[
\begin{aligned}
(|\mathcal B|_{z}^{w})_{j} & < 1 & \forall j,\\
(|\mathcal B|_{z}^{w})_{j} z_{j} (|\mathcal A_{j}|/v)_{i} & < z_{j} (|\mathcal A_{j}|/v)_{i} & \forall j, \ \forall i,\\
(|\mathcal A \mathcal B|_{v}^{w})_{i} \le \sum_{j} (|\mathcal B|_{z}^{w})_{j} z_{j} (|\mathcal A_{j}|/v)_{i} & < \sum_{j} z_{j} (|\mathcal A_{j}|/v)_{i} & \forall i,\\
|\mathcal A \mathcal B|_{v}^{w} & < \sum_{j} z_{j} |\mathcal A_{j}|/v. &
\end{aligned}
\]
Finally, According to Remark \ref{rem:cwmwmn},
\[
|\mathcal A|_{v}^{z} = \sum_{j} z_{j} |\mathcal A_{j}|/v,
\]
which thus concludes the proof.
\endproof
\begin{corollary}
\label{cor:cwmwmn}
Let $\mathcal A$ be a square matrix block-decomposed as
\[
\mathcal A =
\begin{bmatrix}
0 & \cdots & 0 & \mathcal A_{1,m}\\
\vdots & \ddots & \vdots & \vdots\\
0 & \cdots & 0 & \mathcal A_{m-1,m}\\
\mathcal A_{m,1} & \cdots & \mathcal A_{m,m-1} & \mathcal A_{m,m}
\end{bmatrix},
\qquad
m \in \mathbb{N}.
\]
Then, we have
\[
\rho(|\mathcal A|) < 1 \quad \implies \quad \rho(|\mathcal A_{m,m}|+\sum_{i=1}^{m-1}|\mathcal A_{m,i}\mathcal A_{i,m}|) < 1.
\]
\end{corollary}
\proof
Let us consider, for $m = 2$,
\[
\mathcal A =
\begin{bmatrix}
0 & \mathcal A_{1,2}\\
\mathcal A_{2,1} & \mathcal A_{2,2}
\end{bmatrix}.
\]
According to Lemma \ref{lem:pfc},
\[
\rho(|\mathcal A|) < 1 \quad \implies \quad \exists \ w > 0 : \|\mathcal A\|_{\infty}^{w} < 1.
\]
Let such a vector $w$ be decomposed as
\[
w =
\begin{bmatrix}
w_{1} & w_{2}
\end{bmatrix}^{\mathsf T},
\]
and let $u_{1} := (1, 1, ..., 1)^{\mathsf T}$ and $u_{2} := (1, 1, ..., 1)^{\mathsf T}$ be respective size-corresponding vectors. Then, recalling from Remark \ref{rem:cwmwmn} that
\[
\|\mathcal A\|_{\infty}^{w} = \max_{l} |\mathcal A|_{w}^{w}
\]
with $l$ ranging from one to the number of entries of $w$, we also deduce, respectively for the first and the second block-lines of $\mathcal A$, that
\[
\begin{aligned}
\max_{l_{1}} |\mathcal A_{1,2}|_{w_{1}}^{w_{2}} & < 1, &\qquad\qquad \max_{l_{2}} \left(|\mathcal A_{2,1}|_{w_{2}}^{w_{1}} + |\mathcal A_{2,2}|_{w_{2}}^{w_{2}}\right) & < 1,\\
|\mathcal A_{1,2}|_{w_{1}}^{w_{2}} & < u_{1}, &\qquad\qquad |\mathcal A_{2,1}|_{w_{2}}^{w_{1}} + |\mathcal A_{2,2}|_{w_{2}}^{w_{2}} & < u_{2}
\end{aligned}
\]
with $l_{1}$ and $l_{2}$ ranging from one to the number of entries of $w_{1}$ and $w_{2}$, respectively. Lemma \ref{lem:cwmwmn} therefore ensures
\[
\begin{aligned}
|\mathcal A_{2,1} \mathcal A_{1,2}|_{w_{2}}^{w_{2}} & < |\mathcal A_{2,1}|_{w_{2}}^{w_{1}},\\
|\mathcal A_{2,1} \mathcal A_{1,2}|_{w_{2}}^{w_{2}} + |\mathcal A_{2,2}|_{w_{2}}^{w_{2}} & < |\mathcal A_{2,1}|_{w_{2}}^{w_{1}} + |\mathcal A_{2,2}|_{w_{2}}^{w_{2}} < u_{2},\\
\left|(|\mathcal A_{2,1} \mathcal A_{1,2}| + |\mathcal A_{2,2}|)\right|_{w_{2}}^{w_{2}} & < u_{2},\\
\max_{l_{2}} \left|(|\mathcal A_{2,1} \mathcal A_{1,2}| + |\mathcal A_{2,2}|)\right|_{w_{2}}^{w_{2}} & < 1,\\
\left\|(|\mathcal A_{2,1} \mathcal A_{1,2}| + |\mathcal A_{2,2}|)\right\|_{\infty}^{w_{2}} & < 1,\\
\end{aligned}
\]
which concludes the proof for $m = 2$, given that
\[
\rho(|\mathcal A_{2,1} \mathcal A_{1,2}| + |\mathcal A_{2,2}|) \le \left\|(|\mathcal A_{2,1} \mathcal A_{1,2}| + |\mathcal A_{2,2}|)\right\|_{\infty}^{w_{2}}.
\]
The generalization to $m \in \mathbb{N}$ easily follows as we would have, for each block-line $i \in \{1, \ldots, m-1\}$ and the last one, $m$,
\[
|\mathcal A_{i,m}|_{w_{i}}^{w_{m}} < u_{i}, \qquad\qquad \sum_{j = 1}^{m-1}|\mathcal A_{m,j}|_{w_{m}}^{w_{j}} + |\mathcal A_{m,m}|_{w_{m}}^{w_{m}} < u_{m}.
\]
\endproof
Let us now consider
\begin{equation}
\label{eq:ddm:a}
A :=
\begin{bmatrix}
A_{II}^{(1)} & 0 & \cdots & 0 & A_{I\Gamma}^{(1)}\\
0 & A_{II}^{(2)} & \ddots & \vdots & \vdots\\
\vdots & \ddots & \ddots & 0 & A_{I\Gamma}^{(p-1)}\\
0 & \cdots & 0 & A_{II}^{(p)} & A_{I\Gamma}^{(p)}\\
A_{\Gamma I}^{(1)} & \cdots & A_{\Gamma I}^{(p-1)} & A_{\Gamma I}^{(p)} & A_{\Gamma\Gamma}
\end{bmatrix},
\qquad
A_{\Gamma\Gamma} := \sum_{i = 1}^{p} A_{\Gamma\Gamma}^{(i)},
\end{equation}
\begin{equation}
\label{eq:ddm:m}
M :=
\begin{bmatrix}
A_{II}^{(1)} & 0 & \cdots & 0 & 0\\
0 & A_{II}^{(2)} & \ddots & \vdots & \vdots\\
\vdots & \ddots & \ddots & 0 & 0\\
0 & \cdots & 0 & A_{II}^{(p)} & 0\\
0 & \cdots & 0 & 0 & M_{\Gamma\Gamma}
\end{bmatrix}.
\end{equation}
\begin{assumption}
\label{ass:ddm:ai_practical_gen}
Let the matrices $I_{\Gamma\Gamma}^{(i)} - M_{\Gamma\Gamma}^{-1} A_{\Gamma\Gamma}^{(i)}$, $i \in \{1, \ldots, p\}$, be such that
\[
\sum_{i=1}^{p} \left|I_{\Gamma\Gamma}^{(i)} - M_{\Gamma\Gamma}^{-1} A_{\Gamma\Gamma}^{(i)}\right| = \left|\sum_{i=1}^{p} I_{\Gamma\Gamma}^{(i)} - M_{\Gamma\Gamma}^{-1} A_{\Gamma\Gamma}^{(i)}\right|.
\]
\end{assumption}
\begin{theorem}
\label{theo:ddm:ai_practical}
Under Assumption \ref{ass:ddm:ai_practical_gen}, the asynchronous primal Schur method \eqref{eq:ddm:ai} is convergent if
\[
\rho(|I - M^{-1} A|) < 1.
\]
\end{theorem}
\proof
Under Assumption \ref{ass:ddm:ai_practical_gen},
\[
\sum_{i=1}^{p} \left|I_{\Gamma\Gamma}^{(i)} - M_{\Gamma\Gamma}^{-1} A_{\Gamma\Gamma}^{(i)}\right| = \left|I_{\Gamma\Gamma} - M_{\Gamma\Gamma}^{-1} A_{\Gamma\Gamma}^{}\right|,
\]
hence, we have
\[
\begin{aligned}
\sum_{i=1}^{p} \left|I_{\Gamma\Gamma}^{(i)} - M_{\Gamma\Gamma}^{-1} S_{\Gamma\Gamma}^{(i)}\right| & = \sum_{i=1}^{p} \left|I_{\Gamma\Gamma}^{(i)} - M_{\Gamma\Gamma}^{-1} \left(A_{\Gamma\Gamma}^{(i)} - A_{\Gamma I}^{(i)} {A_{II}^{(i)}}^{-1} A_{I\Gamma}^{(i)}\right)\right|\\
& \le \sum_{i=1}^{p} \left|I_{\Gamma\Gamma}^{(i)} - M_{\Gamma\Gamma}^{-1} A_{\Gamma\Gamma}^{(i)}\right| + \left|M_{\Gamma\Gamma}^{-1} A_{\Gamma I}^{(i)} {A_{II}^{(i)}}^{-1} A_{I\Gamma}^{(i)}\right|\\
& = \left|I_{\Gamma\Gamma}^{} - M_{\Gamma\Gamma}^{-1} A_{\Gamma\Gamma}^{}\right| + \sum_{i=1}^{p} \left|M_{\Gamma\Gamma}^{-1} A_{\Gamma I}^{(i)} {A_{II}^{(i)}}^{-1} A_{I\Gamma}^{(i)}\right|.
\end{aligned}
\]
Note, on another side, that
\[
I - M^{-1} A =
\begin{bmatrix}
0 & \cdots & 0 & - {A_{II}^{(1)}}^{-1} A_{I\Gamma}^{(1)}\\
\vdots & \ddots & \vdots & \vdots\\
0 & \cdots & 0 & - {A_{II}^{(p)}}^{-1} A_{I\Gamma}^{(p)}\\
- M_{\Gamma\Gamma}^{-1} A_{\Gamma I}^{(1)} & \cdots & - M_{\Gamma\Gamma}^{-1} A_{\Gamma I}^{(p)} & I_{\Gamma\Gamma}^{} - M_{\Gamma\Gamma}^{-1} A_{\Gamma\Gamma}^{}
\end{bmatrix},
\]
hence, Corollary \ref{cor:cwmwmn} ensures
\[
\rho(|I - M^{-1} A|) < 1 \quad \implies \quad \rho\left(\left|I_{\Gamma\Gamma}^{} - M_{\Gamma\Gamma}^{-1} A_{\Gamma\Gamma}^{}\right| + \sum_{i=1}^{p} \left|M_{\Gamma\Gamma}^{-1} A_{\Gamma I}^{(i)} {A_{II}^{(i)}}^{-1} A_{I\Gamma}^{(i)}\right|\right) < 1.
\]
It follows, according to Lemma \ref{lem:varga2000_theo2.21}, that
\[
\rho\left(\sum_{i=1}^{p} \left|I_{\Gamma\Gamma}^{(i)} - M_{\Gamma\Gamma}^{-1} S_{\Gamma\Gamma}^{(i)}\right|\right) \le \rho\left(\left|I_{\Gamma\Gamma}^{} - M_{\Gamma\Gamma}^{-1} A_{\Gamma\Gamma}^{}\right| + \sum_{i=1}^{p} \left|M_{\Gamma\Gamma}^{-1} A_{\Gamma I}^{(i)} {A_{II}^{(i)}}^{-1} A_{I\Gamma}^{(i)}\right|\right) < 1.
\]
Theorem \ref{theo:ddm:ai} therefore applies, which concludes the proof.
\endproof
\begin{assumption}
\label{ass:ddm:ai_practical}
Let $W_{\Gamma\Gamma}^{(i)} \ge 0$, $i \in \{1, \ldots, p\}$, be matrices of weighting factors such that
\[
A_{\Gamma\Gamma}^{(i)} = W_{\Gamma\Gamma}^{(i)} \circ A_{\Gamma\Gamma}^{},
\qquad
I_{\Gamma\Gamma}^{(i)} = W_{\Gamma\Gamma}^{(i)} \circ I_{\Gamma\Gamma}^{},
\]
where the operator $\circ$ denotes the entry-wise product.
\end{assumption}
\begin{corollary}
\label{cor:ddm:ai_practical}
Under Assumption \ref{ass:ddm:ai_practical}, the asynchronous primal Schur method \eqref{eq:ddm:ai} is convergent if $M_{\Gamma\Gamma}^{}$ is diagonal and
\[
\rho(|I - M^{-1} A|) < 1.
\]
\end{corollary}
\proof
Under Assumption \ref{ass:ddm:ai_practical}, if $M_{\Gamma\Gamma}^{}$ is diagonal, then
\[
\begin{aligned}
I_{\Gamma\Gamma}^{(i)} - M_{\Gamma\Gamma}^{-1} A_{\Gamma\Gamma}^{(i)} & = W_{\Gamma\Gamma}^{(i)} \circ I_{\Gamma\Gamma}^{} - M_{\Gamma\Gamma}^{-1} \left(W_{\Gamma\Gamma}^{(i)} \circ A_{\Gamma\Gamma}^{}\right)\\
& = W_{\Gamma\Gamma}^{(i)} \circ \left(I_{\Gamma\Gamma}^{} - M_{\Gamma\Gamma}^{-1} A_{\Gamma\Gamma}^{}\right).
\end{aligned}
\]
Note that
\[
\sum_{i=1}^{p} A_{\Gamma\Gamma}^{(i)} = A_{\Gamma\Gamma}^{} \quad \implies \quad \sum_{i=1}^{p} W_{\Gamma\Gamma}^{(i)} = (1)_{\Gamma\Gamma}^{},
\]
where $(1)_{\Gamma\Gamma}^{}$ denotes the matrix on $\Gamma$ with all entries set to $1$. It follows that
\[
\sum_{i=1}^{p} \left|I_{\Gamma\Gamma}^{(i)} - M_{\Gamma\Gamma}^{-1} A_{\Gamma\Gamma}^{(i)}\right| = \left|I_{\Gamma\Gamma}^{} - M_{\Gamma\Gamma}^{-1} A_{\Gamma\Gamma}^{}\right|,
\]
hence, Assumption \ref{ass:ddm:ai_practical_gen} is satisfied, which makes Theorem \ref{theo:ddm:ai_practical} applicable.
\endproof
\begin{corollary}
\label{cor:ddm:hsplit}
Under Assumption \ref{ass:ddm:ai_practical_gen}, the asynchronous primal Schur method \eqref{eq:ddm:ai} is convergent if $A$ is an $\mathsf H$-matrix and
\[
\langle M_{\Gamma\Gamma}^{} \rangle - |M_{\Gamma\Gamma}^{} - A_{\Gamma\Gamma}^{}| = \langle A_{\Gamma\Gamma}^{} \rangle.
\]
\end{corollary}
\proof
The matrices $A$ and $M$ being of the forms \eqref{eq:ddm:a} and \eqref{eq:ddm:m}, it follows from Definition \ref{def:compmatrix} of comparison matrices that
\[
\langle M_{\Gamma\Gamma}^{} \rangle - |M_{\Gamma\Gamma}^{} - A_{\Gamma\Gamma}^{}| = \langle A_{\Gamma\Gamma}^{} \rangle
\quad \implies \quad \langle M \rangle - |M - A| = \langle A \rangle.
\]
By Definition \ref{def:hmatrix} of $\mathsf H$-matrices, $\langle A \rangle$ is an $\mathsf M$-matrix,
hence, by Definition \ref{def:hsplit}, $A = M - (M-A)$ is an $\mathsf H$-splitting. Lemma \ref{lem:fs1992} therefore ensures that
\[
\rho(|I - M^{-1} A|) < 1,
\]
hence, Theorem \ref{theo:ddm:ai_practical} applies.
\endproof
\begin{corollary}
\label{cor:ddm:hsplit_diag}
Under Assumption \ref{ass:ddm:ai_practical}, the asynchronous primal Schur method \eqref{eq:ddm:ai} is convergent if $A$ is an $\mathsf H$-matrix, $M_{\Gamma\Gamma}^{}$ is diagonal, and
\[
\langle M_{\Gamma\Gamma}^{} \rangle - |M_{\Gamma\Gamma}^{} - A_{\Gamma\Gamma}^{}| = \langle A_{\Gamma\Gamma}^{} \rangle.
\]
\end{corollary}
\proof
This directly follows from Corollary \ref{cor:ddm:ai_practical} and Corollary \ref{cor:ddm:hsplit}.
\endproof
\begin{notation}
For any square matrix $\mathcal A$, $\operatorname{diag} \mathcal A$ denotes the diagonal matrix obtained by setting all off-diagonal entries of $\mathcal A$ to $0$.
\end{notation}
\begin{remark}
\label{rem:ddm:hsplit}
For any square matrix $\mathcal A$,
\[
\mathcal M = \operatorname{diag} \mathcal A \quad \implies \quad \langle \mathcal M \rangle - |\mathcal M - \mathcal A| = \langle \mathcal A \rangle.
\]
\end{remark}
\begin{remark}
\label{rem:ddm:hsplit2}
For any real square matrix $\mathcal A$ and any real diagonal matrix $\mathcal M$,
\[
\mathcal M \ge \operatorname{diag} \mathcal A \ge 0 \quad \implies \quad \langle \mathcal M \rangle - |\mathcal M - \mathcal A| = \langle \mathcal A \rangle.
\]
\end{remark}

\section{Experimental results}
\label{sec:experiments}

{

\subsection{Implementation aspects}

We followed asynchronous implementation guidelines from \cite{MagGBen2018c} for handling interface data updates, and from \cite{MagGBen2018b,MagGBen2018c,GBenMag2020} for detecting the convergence of the solver. Using the Message Passing Interface (MPI) standard, several message reception requests per neighbor are kept active at the same time, so that several new interface data can be received during computation phase. Communication phase consists of updating the interface buffers using the content of the corresponding communication buffers, then, checking completion of previously initiated requests, upon which new requests are triggered. Similarly, message sending buffers are updated from the content of the newly computed local solution, then message sending requests are triggered if the previous ones have completed. Algorithm~\ref{algo:onelevel_async} gives an overview of the asynchronous primal Schur solver, where the weighting matrix $I^{(i)}$ corresponds to
\[
I^{(i)} := \begin{bmatrix}I & 0\\0 & I^{(i)}_{\Gamma_{l}\Gamma_{l}}\end{bmatrix}.
\]
\begin{algorithm}[htbp]
\caption{Asynchronous primal Schur solver}
\label{algo:onelevel_async}
{
\begin{algorithmic}[1]
\STATE $r^{(i)}$ := $b^{(i)} - A^{(i)} x^{(i)}$
\STATE Synchronize($r^{(i)}_{\Gamma_{l}}$)
\STATE $r^{\mathsf T} r$ := AllReduce(${r^{(i)}}^{\mathsf T}I^{(i)}r^{(i)}$, SUM)
\STATE $\|r\|$ := $\displaystyle\sqrt{r^{\mathsf T} r}$
\STATE $k$ := 0
\STATE $k^{(i)}$ := 0
\STATE convdetect\_phase := 0
\WHILE{$\|r\| > \varepsilon$ \AND $k <$ k\_max}
	\STATE $x^{(i)}_{I}$ := Solve($A_{II}^{(i)}$, $b^{(i)}_{I} - A_{I\Gamma_{l}}^{(i)} x_{\Gamma_{l}}^{(i)}$)
	\STATE $x_{\Gamma_{l}}^{(i)}$ := $I_{\Gamma_{l}\Gamma_{l}}^{(i)} x_{\Gamma_{l}}^{(i)} + R_{\Gamma_{l}\Gamma}^{(i)} M_{\Gamma\Gamma}^{-1} R_{\Gamma\Gamma_{l}}^{(i)}\left(b_{\Gamma_{l}}^{(i)} - A_{\Gamma_{l}I}^{(i)} x_{I}^{(i)} - A_{\Gamma_{l}\Gamma_{l}}^{(i)} x_{\Gamma_{l}}^{(i)}\right)$
	\STATE ASynchronize($x_{\Gamma_{l}}^{(i)}$)
	\STATE $k^{(i)}$ := $k^{(i)} + 1$
	\IF{convdetect\_phase = 0}
			\STATE $r^{(i)}$ := $b^{(i)} - A^{(i)} x^{(i)}$
			\STATE req := ISynchronize($r^{(i)}_{\Gamma_{l}}$)
			\STATE convdetect\_phase := 1
	\ENDIF
	\IF{convdetect\_phase = 1 \AND Test(req)}
			\STATE req := IAllReduce(${r^{(i)}}^{\mathsf T}I^{(i)}r^{(i)}$, $r^{\mathsf T} r$, SUM)
			\STATE convdetect\_phase := 2
	\ENDIF
	\IF{convdetect\_phase = 2 \AND Test(req)}
			\STATE $\|r\|$ := $\displaystyle\sqrt{r^{\mathsf T} r}$
			\STATE convdetect\_phase := 0
			\STATE $k$ := $k + 1$
	\ENDIF
\ENDWHILE
\end{algorithmic}
}
\end{algorithm}
The function call Synchronize($r^{(i)}_{\Gamma_{l}}$) exchanges interface data and performs the update
\[
r^{(i)}_{\Gamma_{l}} := R_{\Gamma_{l}\Gamma}^{(i)} \sum_{j=1}^{p} R_{\Gamma\Gamma_{l}}^{(j)} r_{\Gamma_{l}}^{(j)}.
\]
The function ASynchronize(.) applies the same update without blocking for communication, and managing requests as previously described. The function ISynchronize(.) is non-blocking as well but more similar to MPI non-blocking collective routines such as IAllReduce(.). The update is performed once all of the interface data is available, then the returned request object is marked as having completed. Such non-blocking synchronization is derived from snapshot-based convergence detection \cite{MagGBen2018b,MagGBen2018c}, while being able to rely on such residual error computed without the exact snapshot-based protocol is due to \cite{GBenMag2020}.

}

\subsection{Problems and settings}

\begin{figure}[htbp]
  \centering
  \includegraphics[scale=0.18]{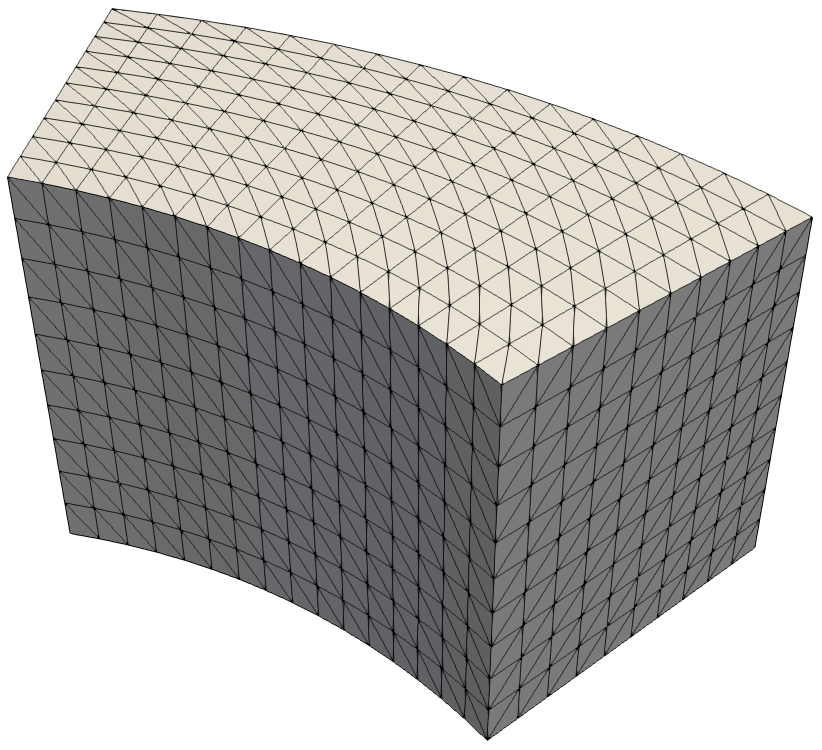}
  \qquad
  \qquad
  \includegraphics[scale=0.18]{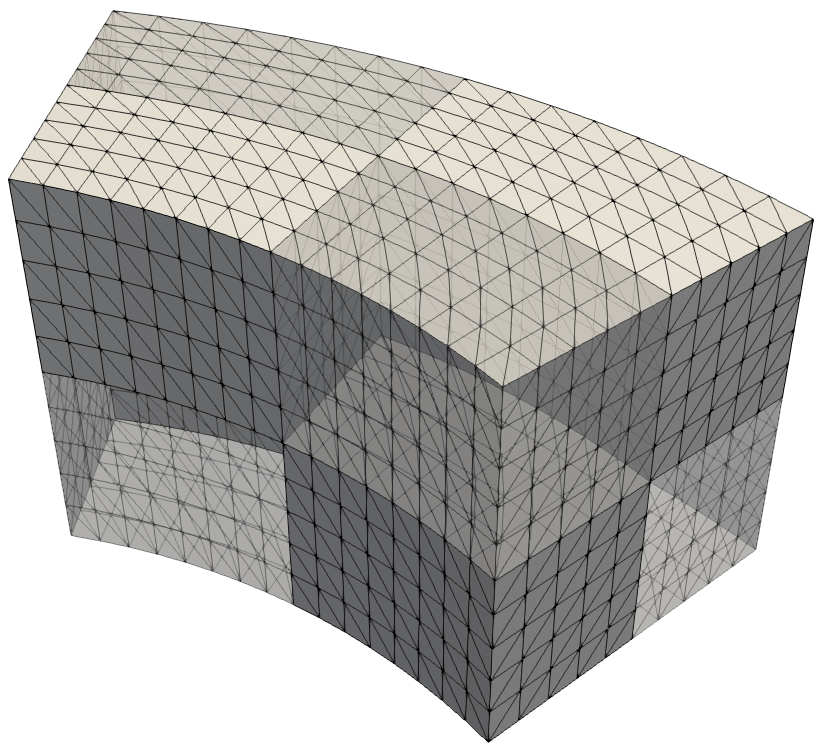}
  \caption{(left) Experimental domain: right helicoid of length $\pi/4$, (right) 3D partitioning scheme.}
  \label{fig:res:dom}
\end{figure}

\begin{figure}[htbp]
  \centering
  \includegraphics[scale=0.18]{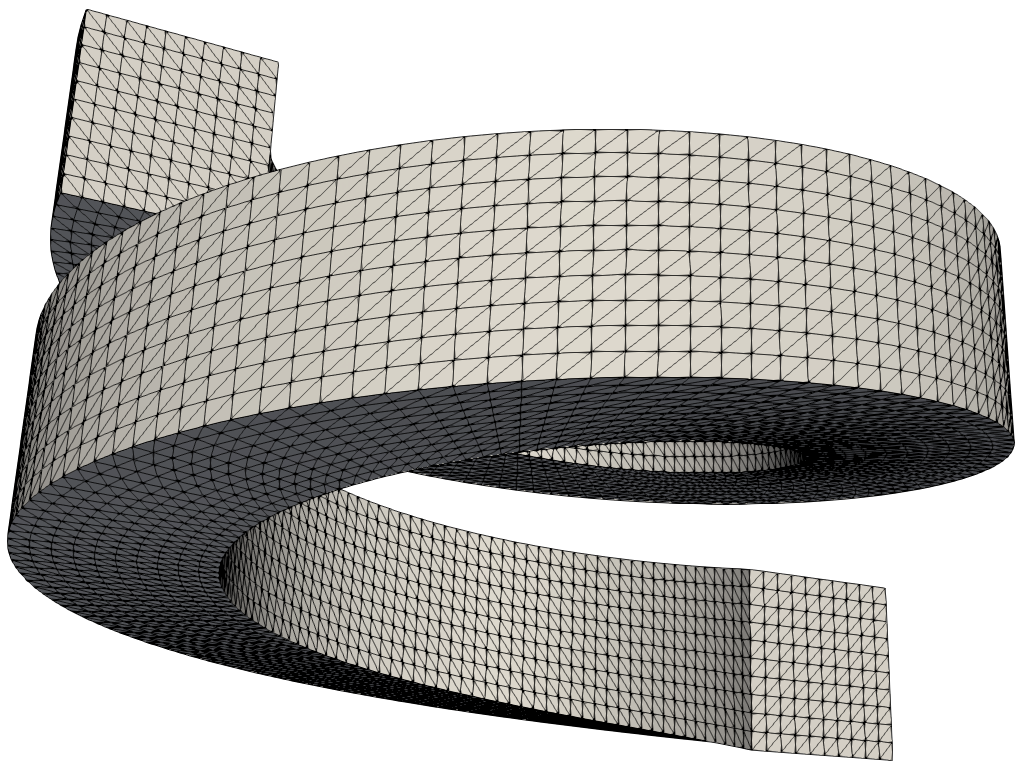}
  \caption{Right helicoid of length $3 \pi$.}
  \label{fig:res:helic}
\end{figure}

As numerical experiment, we consider a right helicoid domain shown in Figure \ref{fig:res:dom}, left,
\graphicspath{{./}} %Setting the graphicspath
which dimensions consist of its height $h = 1$~m in the $z$-axis direction, its width $w = 1$~m in the $x$-axis direction and its length $l = \pi/4$ denoting the total rotation covered. Figure \ref{fig:res:helic}, for instance, shows a helicoid of length $3 \pi$.
A last characteristic is the raising factor of the helicoid while rotating, which is set here to $0.25$. See, e.g., \cite{KrivoRyn2017} for geometrical equations fully describing different types of helicoid. The partitioning of the geometry is made along each of its dimensions, as shown in Figure \ref{fig:res:dom}, right.

\begin{figure}[htbp]
  \centering
  \includegraphics[scale=0.18]{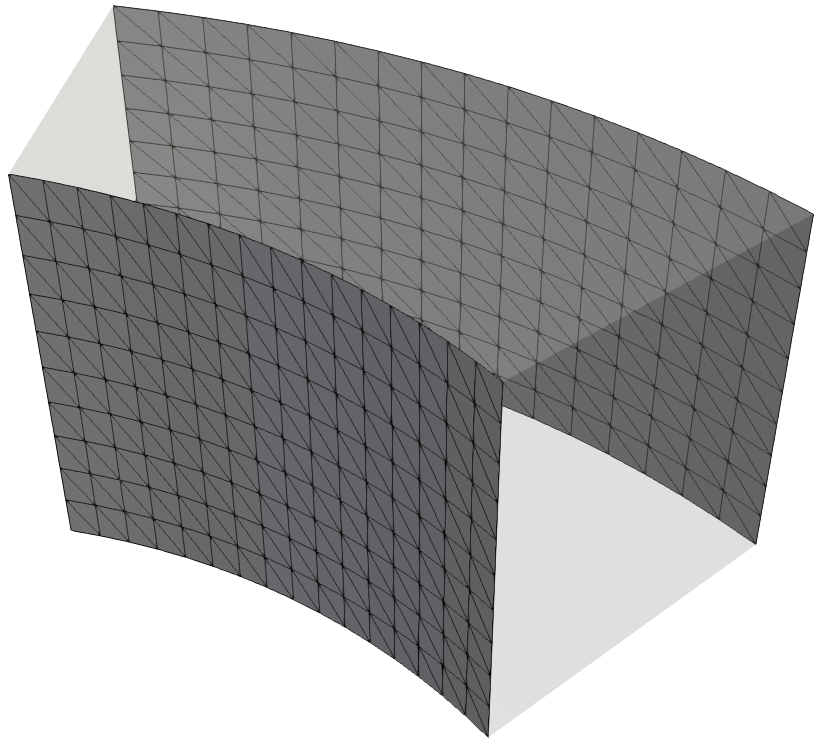}
  \qquad
  \qquad
  \includegraphics[scale=0.18]{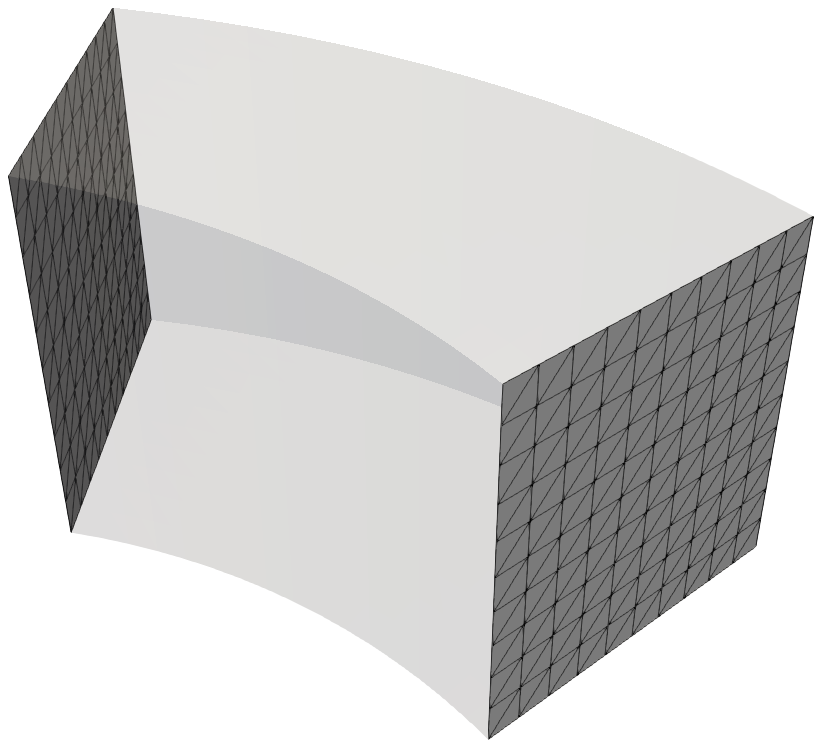}
  \caption{Boundary prescription: (left) inner and outer sides, (right) starting and ending sides.}
  \label{fig:res:bnd}
\end{figure}

Both Poisson's problem,
\[
- \Delta u = g,
\]
and linear elasticity,
\[
- \operatorname{div} \sigma(u) = g,
\]
are investigated with a uniform volume source (or load) and prescribed solution, $u = 0$, on either the inner and outer sides of the helicoid, or its starting and ending sides (see Figure \ref{fig:res:bnd}). In the linear elasticity case, a concrete material is considered with Young's modulus $E = 32.5$~GPa and Poisson's ratio $\nu = 0.18$.

In each subdomain, stiffness matrix and right-hand-side vector are generated by P1 finite-element approximation. According to Remark \ref{rem:ddm:hsplit2}, splitting matrices
\[
M_{\Gamma\Gamma}^{} = \alpha \operatorname{diag} A_{\Gamma\Gamma}^{}, \qquad \alpha \ge 1,
\]
are applied with $\alpha := 1$ for the Poisson's problem, and $\alpha := 1.4$ for the linear elasticity problem.

The experimental environment consists of an SGI ICE X cluster with 30 computing nodes, where each node contains two 12-cores Intel Haswell Xeon CPUs at 2.30~GHz, and 58~GB RAM allocated, leading to a maximum of 720 processor cores available. The nodes are interconnected through an FDR Infiniband network (56~Gb/s), and communication is handled by the SGI-MPT middleware which implements the MPI standard. Each processor core is mapped with one MPI process, which is mapped with one subdomain.

In the sequel, we present execution times $t$ in seconds, numbers of iterations $k$ and final global residual errors $\|A x - b\|$, where $\|.\|$ denotes the Euclidean norm, for a { residual} threshold set to 1e-6 { as stopping criterion}. Comparison is made between conjugate gradient iterations (CG-Schur) and our relaxation approach (async-Schur), both combined with $\mathsf{LU}$-factorization of the local matrices $A_{II}^{(i)}$. In case of asynchronous processing, each process $i$ performs a proper number $k^{(i)}$ of iterations, then $k_{\max} := \max_{i} k^{(i)}$ will be considered here. The final residual error is calculated synchronously after convergence is detected. Combining the snapshot-based and protocol-free approaches from \cite{MagGBen2018b, MagGBen2018c} and \cite{GBenMag2020}, intermediate global residual errors could be evaluated conjointly to the asynchronous iterations without slowing down the processes, leading to consistently quick and accurate convergence detection.

\subsection{Performance assessment}

The first observations are made on the Poisson's problem with Dirichlet boundary conditions set on the inner and outer sides of the helicoid (see Figure \ref{fig:res:bnd}, left). Table \ref{tab:res:poisson_io_120} shows a limit of the problem size $n$ (number of degrees of freedom) beyond which the performance of the CG-Schur method is strongly impacted, contrarily to our async-Schur approach which features a relatively very low increase of its execution time when switching from subdomain's average size $n^{(i)} = 1139$ to $n^{(i)} = 1572$. From there, async-Schur stays relatively close to CG-Schur, and even outperforms it for $n^{(i)} = 1572$.
\begin{table}[htbp]
\caption{Performance for $- \Delta u = g$, $u = 0$ on the helicoid's inner and outer sides, 120 cores (5 nodes).}
\label{tab:res:poisson_io_120}
\centering
{\footnotesize
\begin{tabular}{ccc}
\hline\noalign{\smallskip}
& CG-Schur & Async-Schur\\
\noalign{\smallskip}\hline\noalign{\smallskip}
\begin{tabular}{rr}
$n$ & $n^{(i)}$\\
\noalign{\smallskip}\hline\noalign{\smallskip}
72576 & 793\\
107584 & \textbf{1139}\\
152352 & \textbf{1572}\\
208080 & 2104\\
275968 & 2744\\
357216 & 3501\\
\end{tabular}
&
\begin{tabular}{rrr}
$t$ (sec) & $k$ & $\|A x - b\|$\\
\noalign{\smallskip}\hline\noalign{\smallskip}
0.2 & 135 & 8.99E-07\\
\textbf{0.3} & 146 & 9.26E-07\\
\textbf{21} & 152 & 9.38E-07\\
41 & 160 & 9.35E-07\\
55 & 171 & 9.67E-07\\
67 & 174 & 8.69E-07\\
\end{tabular}
&
\begin{tabular}{rrrr}
$t$ (sec) & $k_{\max}$ & $\|A x - b\|$ & $\mathbf{t / t_{\mathtt{\textbf{CG}}}}$\\
\noalign{\smallskip}\hline\noalign{\smallskip}
5 & 12101 & 9.24E-07 & \textbf{22}\\
\textbf{6} & 12079 & 9.76E-07 & \textbf{23}\\
\textbf{19} & 15942 & 8.51E-07 & \textbf{0.89}\\
70 & 28747 & 9.51E-07 & \textbf{1.69}\\
129 & 16846 & 9.67E-07 & \textbf{2.34}\\
265 & 17502 & 9.77E-07 & \textbf{3.97}\\
\end{tabular}\\
\noalign{\smallskip}\hline
\end{tabular}
}
\end{table}

{

For average subdomain's sizes $n^{(i)} > 1139$, Figure \ref{fig:res:scaling} gives a performance overview when increasing the number of subdomains but keeping a constant problem size $n$ (weak scaling). This also corresponds to a decreasing ratio $n^{(i)}/n$.
\begin{figure}[htbp]
  \centering
  \includegraphics[scale=0.37]{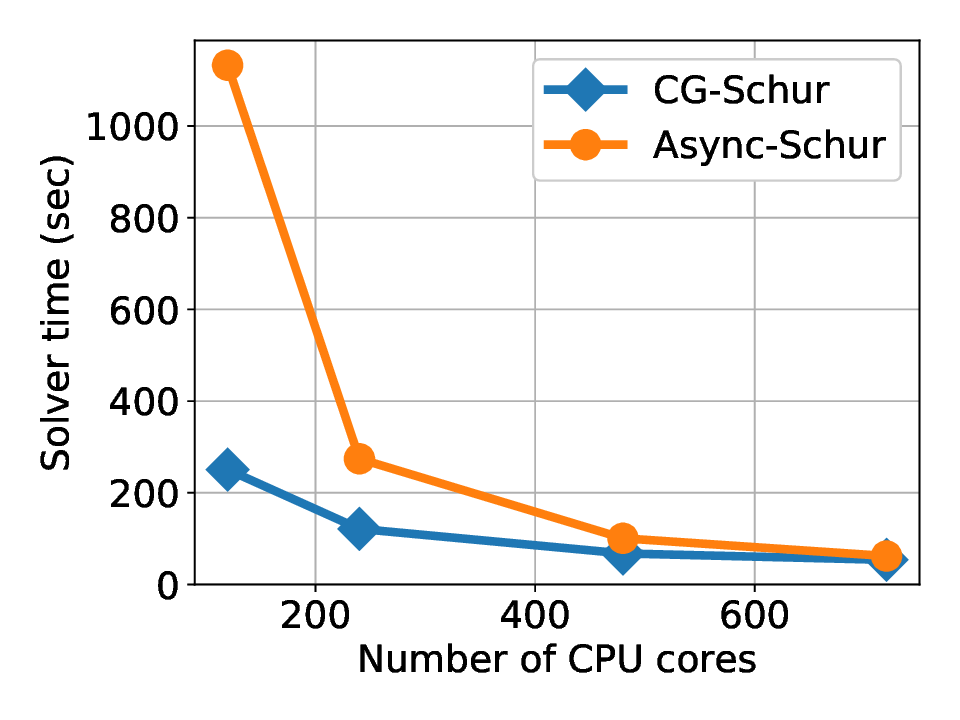}
  \caption{ Weak scaling performances for $- \Delta u = g$, $u = 0$ on the  helicoid's inner and outer sides, problem size $n = 833247$.}
  \label{fig:res:scaling}
\end{figure}
Considering, however, observations illustrated in Table \ref{tab:res:poisson_io_120}, it remains very likely that the performance of the async-Schur method depends more on $n^{(i)}$ than on $n^{(i)}/n$.

}

Such possible close performances of the CG-Schur and the async-Schur make it interesting to assess the potential impact of failures, given that their probability increases with the number of computing nodes. According to the speedup values, we simulate 5 temporary failures successively occurring on different computing nodes, each time after approximately 90\% of the fault-free CG execution time (see Table \ref{tab:res:elast_ft_sim}).
\begin{table}[htbp]
\caption{Failure occurrence time (in seconds) during solver execution, for $p \in \{120, 240, 480, 720\}$, $- \operatorname{div} \sigma(u) = g$, $u = 0$ on the helicoid's starting and ending sides.}
\label{tab:res:elast_ft_sim}
\centering
{\footnotesize
\begin{tabular}{ccc}
\hline\noalign{\smallskip}
& CG-Schur & Async-Schur\\
\noalign{\smallskip}\hline\noalign{\smallskip}
\begin{tabular}{c}
Failed node\\
\noalign{\smallskip}\hline\noalign{\smallskip}
0\\
1\\
2\\
3\\
4\\
\end{tabular}
&
\begin{tabular}{rrrr}
120 & 240 & 480 & 720\\
\noalign{\smallskip}\hline\noalign{\smallskip}
56 & 76 & 116 & 138\\
113 & 152 & 232 & 275\\
169 & 228 & 346 & 412\\
225 & 303 & 462 & 550\\
280 & 379 & 576 & 687\\
\end{tabular}
&
\begin{tabular}{rrrr}
120 & 240 & 480 & 720\\
\noalign{\smallskip}\hline\noalign{\smallskip}
55 & 73 & 111 & 132\\
109 & 145 & 220 & 266\\
164 & 216 & 330 & 403\\
218 & 288 & 441 & 546\\
272 & 360 & 552 & 691\\
\end{tabular}\\
\noalign{\smallskip}\hline
\end{tabular}
}
\end{table}
For the 24 subdomains being processed on a node, a failure (and recovery) is simulated by resetting $x_{\Gamma_{l}}^{(i),k}$ and communication buffers to their initial value. Factorization of $A_{II}^{(i)}$ is assumed to have been backed up, hence, is not reprocessed. Then, for the CG, the global iterative procedure has to synchronize and restart, while nothing needs to be particularly done in case of asynchronous iterations. Table \ref{tab:res:elast_ft} confirms the small impact on the asynchronous solver, while restarting the CG iterations do not seem to however benefit from the initial guess partially very close to the solution (from the fault-free subdomains), no matter what is the proportion of this advanced part in the global interface vector.
\begin{table}[htbp]
\caption{Performance with 5 node failures, for $- \operatorname{div} \sigma(u) = g$, $u = 0$ on the helicoid's starting and ending sides.}
\label{tab:res:elast_ft}
\centering
{\footnotesize
\begin{tabular}{ccc}
\hline\noalign{\smallskip}
& CG-Schur & Async-Schur\\
\noalign{\smallskip}\hline\noalign{\smallskip}
\begin{tabular}{crr}
$p$ \textbar\ \#Nodes & Failed & $n^{(i)}$\\
& nodes &\\
\noalign{\smallskip}\hline\noalign{\smallskip}
120 \textbar\ 05 & 100\% & 1898\\
240 \textbar\ 10 & 50\% & 1973\\
480 \textbar\ 20 & 25\% & 2001\\
720 \textbar\ 30 & 17\% & 2053\\
\end{tabular}
&
\begin{tabular}{ccr}
$t$ & $k$ & $\|A x - b\|$\\
(sec) &&\\
\noalign{\smallskip}\hline\noalign{\smallskip}
347 & 1901 & 9.40E-07\\
471 & 2451 & 9.16E-07\\
710 & 3136 & 9.26E-07\\
847 & 3613 & 9.97E-07\\
\end{tabular}
&
\begin{tabular}{ccr}
$t$ & $k_{\max}$ & $\|A x - b\|$\\
(sec) &&\\
\noalign{\smallskip}\hline\noalign{\smallskip}
344 & 373817 & 9.53E-07\\
456 & 482740 & 9.77E-07\\
689 & 744514 & 9.05E-07\\
886 & 984099 & 9.57E-07\\
\end{tabular}\\
\noalign{\smallskip}\hline
\end{tabular}
}
\end{table}

{

It is a big challenge to make relaxation methods competitive with the CG method. While still not superior, asynchronous relaxation in the domain decomposition framework constitutes a huge step towards it, compared to its synchronous counterpart. We see from the third line of Table \ref{tab:res:poisson_io_120} that cases can be found where the CG is very slightly outperformed. While not reported here, a ratio $t_{\mathtt{async}} / t_{\mathtt{CG}}$ of 0.94 was also observed with $p = 240$ and $n^{(i)} = 1449$. However, in all of our other experienced configurations ($p \ge 480$ or other prescribed boundaries), best ratios were found between 1.1 and 2. As commonly claimed though, Table \ref{tab:res:elast_ft} confirms the possibility to switch to the async-Schur method in case of high probability or frequency of node failure.

}

\section{Conclusions}
\label{sec:conclusions}

The asynchronous iterations theory was developed as a generalization of classical relaxation schemes, and, so far, domain decomposition methods are investigated in their asynchronous version only when they already exhibit an inherent relaxation framework. The classical primal Schur method does not directly provide such a feature, which has made it require, first, to model an applicable relaxation scheme within its framework. Second, for the usual context where the local Schur complement matrices are not actually constructed, suitable matrix splittings have been provided, only based on matrices explicitly generated. While this designed asynchronous primal Schur method constitutes an unusual application of the asynchronous iterations theory, its convergence has been established under classical conditions. Experimental results have shown that, depending on the size of the subdomains, the relaxation-based asynchronous iterative solver can be competitive with a conjugate gradient one, especially in case of hardware failures occurring in iterations loops with no backup of the updated vectors. Being based on a quite unusual general asynchronous fixed-point theory due to Baudet, the present work possibly opens new developments in the asynchronous iterative computing field.

\section*{Acknowledgement}

The research has been supported by the ``RUDN University Program 5-100'', the French national program LEFE/INSU, the project ADOM (M\'ethodes de d\'ecomposition de domaine asynchrones) of the French National Research Agency (ANR), and using HPC resources from the ``M\'esocentre'' computing center of CentraleSup\'elec and \'Ecole Normale Sup\'erieure Paris-Saclay supported by CNRS and R\'egion \^Ile-de-France.

\bibliography{ref}

\begin{thebibliography}{10}

\bibitem{Amdahl1967}
G.~M. Amdahl.
\newblock Validity of the single processor approach to achieving large scale
  computing capabilities.
\newblock In {\em Proceedings of the April 18-20, 1967, Spring Joint Computer
  Conference, Atlantic City, New Jersey, USA}, AFIPS '67 (Spring), pages
  483--485, New York, NY, USA, 1967. ACM.

\bibitem{AxelKolo1994}
O.~Axelsson and L.~Kolotilina.
\newblock Diagonally compensated reduction and related preconditioning methods.
\newblock {\em Numer. Linear Algebra Appl.}, 1(2):155--177, 1994.

\bibitem{BahiEtAl1996}
J.~Bahi, E.~Griepentrog, and J.~C. Miellou.
\newblock Parallel treatment of a class of differential-algebraic systems.
\newblock {\em SIAM J. Numer. Anal.}, 33(5):1969--1980, 1996.

\bibitem{Baudet1978}
G.~M. Baudet.
\newblock Asynchronous iterative methods for multiprocessors.
\newblock {\em J. ACM}, 25(2):226--244, 1978.

\bibitem{BertTsit1989}
D.~P. Bertsekas and J.~N. Tsitsiklis.
\newblock {\em Parallel and Distributed Computation: Numerical Methods}.
\newblock Prentice-Hall, Inc., Upper Saddle River, NJ, USA, 1989.

\bibitem{BertTsit1991}
D.~P. Bertsekas and J.~N. Tsitsiklis.
\newblock Some aspects of parallel and distributed iterative algorithms -- {A}
  survey.
\newblock {\em Automatica J. IFAC}, 27(1):3--21, 1991.

\bibitem{CastelEtAl1998}
M.~J. Castel, V.~Migall\'on, and J.~Penad\'es.
\newblock Convergence of non-stationary parallel multisplitting methods for
  hermitian positive definite matrices.
\newblock {\em Math. Comp.}, 67(221):209--220, 1998.

\bibitem{ChazMir1969}
D.~Chazan and W.~Miranker.
\newblock Chaotic relaxation.
\newblock {\em Linear Algebra Appl.}, 2(2):199--222, 1969.

\bibitem{Fan1958}
K.~Fan.
\newblock Topological proofs for certain theorems on matrices with non-negative
  elements.
\newblock {\em Monatshefte f\"ur Mathematik}, 62:219--237, 1958.

\bibitem{Fan1960}
K.~Fan.
\newblock Note on {$M$-matrices}.
\newblock {\em Q. J. Math.}, 11(1):43--49, 1960.

\bibitem{FromEtAl1997}
A.~Frommer, H.~Schwandt, and D.~B. Szyld.
\newblock Asynchronous weighted additive {Schwarz} methods.
\newblock {\em Electron. Trans. Numer. Anal.}, 5:48--61, 1997.

\bibitem{FromSzyld1992}
A.~Frommer and D.~B. Szyld.
\newblock H-splittings and two-stage iterative methods.
\newblock {\em Numer. Math.}, 63(1):345--356, 1992.

\bibitem{GBenMag2020}
G.~Gbikpi-Benissan and F.~Magoul\`es.
\newblock Protocol-free asynchronous iterations termination.
\newblock {\em Adv. Eng. Softw.}, 146:102827, 2020.

\bibitem{KrivoRyn2017}
{Krivoshapko, S. N.} and {Rynkovskaya, M.}
\newblock Five types of ruled helical surfaces for helical conveyers, support
  anchors and screws.
\newblock {\em MATEC Web Conf.}, 95:06002, 2017.

\bibitem{MagGBen2018}
F.~Magoul\`es and G.~Gbikpi-Benissan.
\newblock Asynchronous parareal time discretization for partial differential
  equations.
\newblock {\em SIAM J. Sci. Comput.}, 40(6):C704--C725, 2018.

\bibitem{MagGBen2018b}
F.~Magoul\`es and G.~Gbikpi-Benissan.
\newblock Distributed convergence detection based on global residual error
  under asynchronous iterations.
\newblock {\em IEEE Trans. Parallel Distrib. Syst.}, 29(4):819--829, 2018.

\bibitem{MagGBen2018c}
F.~Magoul\`es and G.~Gbikpi-Benissan.
\newblock {JACK2}: An {MPI}-based communication library with non-blocking
  synchronization for asynchronous iterations.
\newblock {\em Adv. Eng. Softw.}, 119:116--133, 2018.

\bibitem{MagoulesEtAl2004}
F.~Magoul\`es, F.-X. Roux, and S.~Salmon.
\newblock Optimal discrete transmission conditions for a nonoverlapping domain
  decomposition method for the {Helmholtz} equation.
\newblock {\em SIAM J. Sci. Comput.}, 25(5):1497--1515, 2004.

\bibitem{MagEtAl2017}
F.~Magoul{\`e}s, D.~B. Szyld, and C.~Venet.
\newblock Asynchronous optimized {S}chwarz methods with and without overlap.
\newblock {\em Numer. Math.}, 137(1):199--227, 2017.

\bibitem{MagVen2018}
F.~Magoul\`es and C.~Venet.
\newblock Asynchronous iterative sub-structuring methods.
\newblock {\em Math. Comput. Simulation}, 145(Supplement C):34--49, 2018.

\bibitem{Schechter1959}
S.~Schechter.
\newblock Relaxation methods for linear equations.
\newblock {\em Comm. Pure Appl. Math.}, 12(2):313--335, 1959.

\bibitem{SpitEtAl2001}
P.~Spit\'eri, J.-C. Miellou, and D.~El~Baz.
\newblock Asynchronous {Schwarz} alternating method with flexible communication
  for the obstacle problem.
\newblock {\em R\'es. Syst. R\'epartis Calculateurs Parall\`eles},
  13(1):47--66, 2001.

\bibitem{SpitEtAl2003}
P.~Spit\'eri, J.-C. Miellou, and D.~El~Baz.
\newblock Parallel asynchronous {S}chwarz and multisplitting methods for a
  nonlinear diffusion problem.
\newblock {\em Numer. Algorithms}, 33(1):461--474, 2003.

\bibitem{Varga2000}
R.~S. Varga.
\newblock {\em Matrix Iterative Analysis}, volume~27 of {\em Springer Series in
  Computational Mathematics}.
\newblock Springer-Verlag Berlin Heidelberg, 2000.

\bibitem{YunKim2004}
J.~H. Yun and S.~W. Kim.
\newblock Convergence of two-stage iterative methods using incomplete
  factorization.
\newblock {\em J. Comput. Appl. Math.}, 166(2):565--580, 2004.

\end{thebibliography}
\bibliographystyle{abbrv}

\end{document}